\documentclass{hajlasz}

\usepackage{amssymb}
\newtheorem{theorem}{Theorem}[section]
\newtheorem{lemma}[theorem]{Lemma}

\theoremstyle{definition}

\theoremstyle{remark}

\numberwithin{equation}{section}

\def\id{{\rm id\, }}
\def\eps{\varepsilon}
\def\Om{\Omega}

\def\vi{\varphi}

\def\M{{\mathcal M}}
\def\H{{\mathcal H}}
\def\charfn_#1{{\raise1.2pt\hbox{$\chi_{\kern-1pt\lower3pt\hbox{{$\scriptstyle#1$}}}$}}}
\def\charfun_#1{{\raise1.2pt\hbox{$\chi_{\kern-1pt\lower6pt\hbox{{$\scriptstyle#1$}}}$}}}
\def\mvint_#1{\mathchoice
          {\mathop{\vrule width 6pt height 3 pt depth -2.5pt
                  \kern -8pt \intop}\nolimits_{\kern -3pt #1}}%
          {\mathop{\vrule width 5pt height 3 pt depth -2.6pt
                  \kern -6pt \intop}\nolimits_{#1}}%
          {\mathop{\vrule width 5pt height 3 pt depth -2.6pt
                  \kern -6pt \intop}\nolimits_{#1}}%
          {\mathop{\vrule width 5pt height 3 pt depth -2.6pt
                  \kern -6pt \intop}\nolimits_{#1}}}

\newfont{\pogrubianemat}{msbm10}%
\newfont{\malepogrubianemat}{msbm7}%

\def\bbbr{\mbox{\pogrubianemat R}}

plus 6pt minus 12pt
\def\dist{{\rm{dist}\,}}              

\def\lip{{\rm{Lip}\,}}                
\def\supp{{\rm{supp}\,}}              

\def\dimh{\hbox{${\rm dim}_{\lower1.5pt\hbox{{$\scriptstyle H$}}
                         }\kern.7pt$}
        }
\def\rank{{{\rm rank}\kern.7pt}}

\newcounter{cond}

\begin{document}

\sloppy
\frenchspacing

\title{Sobolev mappings: Lipschitz density is not a
bi-Lipschitz invariant of the target}

\author{Piotr Haj\l{}asz}
\address{Department of Mathematics, University of
Pittsburgh, 301 Thackeray Hall, Pittsburgh, PA 15260, USA}
\email{hajlasz@pitt.edu}

\thanks{This work was supported by the NSF grant DMS-0500966.}

\date{}


\maketitle

\pagestyle{myheadings}
\markboth{\normalsize Piotr Haj\l{}asz}{\normalsize Lipschitz density is not a bi-Lipschitz invariant of the target}

\begin{abstract}
We study a question of density of Lipschitz mappings in the
Sobolev class of mappings from a closed manifold into a singular
space. The main result of the
paper, Theorem~\ref{main}, shows that if we change the metric in
the target space to a bi-Lipschitz equivalent one, than the property of
the density of Lipschitz mappings may be lost.
Other main results in the paper are
Theorems~\ref{one}, \ref{dens}, \ref{subaru}, \ref{smooth}.
\end{abstract}

\section{Introduction}

Sobolev spaces have been playing a major role in the study of
nonlinear partial differential equations
and calculus of variations for more than seven
decades now. Variational approach to mappings between manifolds
leads in a natural way to the notion of Sobolev mappings between
manifolds $W^{1,p}(M,N)$.
Here the symbol $W^{1,p}$ stands for the Sobolev
space of $L^p$ integrable
mappings with $L^p$ integrable first order derivatives.
This class of mappings was
considered for the first time as early as in 1948 by Morrey
\cite{morrey} in his celebrated paper on regularity of minimizing harmonic
mappings from two dimensional surfaces into manifolds.
While it is a well known theorem of Meyers and Serrin
that smooth functions are dense in the Sobolev space $W^{1,p}(M)$,
where $M$ is a Riemannian manifold,
similar
question, raised by Eells and Lemaire, \cite{eellsl},
about density of smooth mappings in the class of Sobolev
mappings between manifolds is no longer obvious.
Actually Schoen and Uhlenbeck \cite{schoenu1}, \cite{schoenu2}
were the first to observe that in general the answer to this question is in the
negative. Namely they demonstrated
that the radial projection mapping
$x/|x|\in W^{1,p}(B^n,S^{n-1})$, $1\leq p<n$,  cannot be
approximated by $C^{\infty}(B^n,S^{n-1})$ mappings when $n-1\leq p<n$. In the
same paper they proved that $C^{\infty}(M,N)$ is a dense subset of
$W^{1,p}(M,N)$ when $p\geq\dim M$ and $M$, $N$ are compact smooth manifolds,
$\partial N=\emptyset$.
Later, mainly due to the work of Bethuel, \cite{bethuel}, it turned
out that the condition for the density of
smooth mappings can be formulated in terms of algebraic topology
describing
the topological structure of the manifolds $M$ and $N$.
Finally a necessary and sufficient condition for the density of
smooth mappings in the space of Sobolev mappings $W^{1,p}(M,N)$
was discovered by Hang and Lin \cite{hangl1}, \cite{hangl2}.
This result corrects an earlier statement of Bethuel \cite{bethuel}
about a necessary and sufficient condition for the density.
Other papers on this topic include
\cite{bethuel1}, \cite{bethuelcdh},
\cite{bethuelz},
\cite{brezisl}, \cite{brezislmn}, \cite{brezisn1}, \cite{brezisn2}
\cite{corong}, \cite{demengel}, \cite{escobedo}, \cite{grecoiss},
\cite{hajlasz1}, \cite{hajlasz2}, \cite{HIMO},
\cite{hangl1}, \cite{hangl2}, \cite{helein1},
\cite{isobe},
\cite{pakzadr}, \cite{white}, \cite{white2}.

The theory of Sobolev mappings between manifolds has been extended to the case
of Sobolev mappings with values into metric spaces.
The first papers on this subject include the work of Ambrosio,
\cite{ambrosio}, on limits of classical variational problems
and the work of Gromov and Schoen, \cite{gromovs}, on
Sobolev mappings into the Bruhat-Tits buildings,
with applications to rigidity questions for discrete groups.
Later the theory of Sobolev mappings with values into metric
spaces was developed in a more elaborated form by Korevaar and Schoen
\cite{korevaars} in their
approach to the theory of harmonic mappings into
Alexandrov spaces of non-positive curvature.
Other papers on Sobolev mappings from a manifold into
a metric space include \cite{capognal}, \cite{eellsf},
\cite{jost1}, \cite{jost2}, \cite{jost3}, \cite{jostz},
\cite{reshetnyak}, \cite{serbinowski}.
Finally analysis on metric spaces,
the theory of Carnot--Carath\'eodory spaces and
the theory of quasiconformal mappings between metric spaces lead to
the theory of Sobolev mappings between
metric spaces, \cite{heinonenk}, \cite{heinonenkst},  \cite{troyanov}.

In the case of Sobolev mappings into metric spaces Heinonen,
Koskela, Shanumgalingam and Tyson, \cite[Remark~6.9]{heinonenkst},
ask for the conditions that would guarantee the density of
Lipschitz mappings.
Our paper is devoted to study this problem.
From the point of view of geometric analysis one is willing to
identify spaces that are bi-Lipschitz equivalent.
The main result of the
paper, Theorem~\ref{main}, shows that if we change the metric in
the target space to a bi-Lipschitz equivalent metric, than the property of
the density of Lipschitz mappings may be lost. That shows a
major difficulty in attempts to develop the theory in the case of
metric space valued maps. Other main results in the paper are
Theorems~\ref{one}, \ref{dens}, \ref{subaru}, \ref{smooth}. All
the main results but Theorem~\ref{dens} are counterexamples.

Let $M$ and $N$ be two smooth compact manifolds, $N$ without boundary. Assume
that $N$ is smoothly embedded into
the Euclidean space $\bbbr^\nu$.
Then for $1\leq p<\infty$ we define the class of
Sobolev mappings between manifolds $M$ and $N$
as
\begin{equation}
\label{betweenman}
W^{1,p}(M,N) =
\{ u\in W^{1,p}(M,\bbbr^\nu):\, u(x)\in N\ \text{a.e.} \}.
\end{equation}
The space $W^{1,p}(M,N)$ is equipped with the metric of the Sobolev norm
$W^{1,p}$ i.e. $\rho(u,v)=\Vert u-v\Vert_{1,p}$ for $u,v\in W^{1,p}(M,N)$.

By the theorem of Meyers and Serrin every Sobolev mapping
$u\in W^{1,p}(M,N)$ can be approximated by smooth mappings
$C^\infty(M,\bbbr^\nu)$, but it is not always possible to choose
the approximating sequence to have values in $N$.

Clearly the space $W^{1,p}(M,N)$ depends on the embedding of $N$ into
the Euclidean
space. However, if $i_1:N\to\bbbr^{\nu_1}$ and $i_2:N\to\bbbr^{\nu_2}$ are two
different smooth embeddings  then the mapping
$$
\Phi:W^{1,p}(M,i_1(N))\to W^{1,p}(M,i_2(N)),
\qquad
\Phi(u) = i_2\circ i^{-1}_{1}\circ u
$$
is a homeomorphism of the spaces $W^{1,p}(M,i_1(N))$ and
$W^{1,p}(M,i_2(N))$.
In particular the answer to the question about density of smooth mappings in
$W^{1,p}(M,N)$ does not depend on the smooth embedding of $N$.

There are various possibilities to define metric space valued
Sobolev mappings, see e.g. \cite{heinonenkst}, \cite{korevaars}, \cite{reshetnyak}.
As was proved
by Heinonen, Koskela, Shanmugalingam and Tyson
\cite{heinonenkst}, in a
reasonable generality different approaches are equivalent.
The results in \cite{heinonenkst} concern
Sobolev mappings between metric spaces, but in our paper we will
consider the case of mappings from a closed manifold to metric
spaces $W^{1,p}(M,X)$ only.
One of the definitions of such mappings mimics  the
definition (\ref{betweenman}). First one defines Sobolev mappings
from $M$ into Banach spaces. Then one considers an isometric
embedding of $X$ into a Banach space $V$
(any metric space $X$ admits an isometric embedding into
$\ell^\infty(X)$)
and defines
$W^{1,p}(M,X)$ as the set of all mappings in $W^{1,p}(M,V)$ such
that values of the mapping belong to $X\subset V$. For a precise
definition and proof of the equivalence with other approaches, see
\cite{heinonenkst}.
The space $W^{1,p}(M,X)$ inherits the metric from the Banach
structure of the space $W^{1,p}(M,V)$. Then under reasonable
assumptions
every Sobolev mapping
$u\in W^{1,p}(M,X)$ can be approximated by Lipschitz mappings
$\lip(M,V)$,
\cite[Theorem~6.7 and Remark~6.6]{heinonenkst}.
The authors ask:
{\em It is an interesting problem to determine when one can choose
the Lipschitz approximation to have values in the target $X$.
[\ldots] For instance, one can ask to what extent Bethuel's
results have analogs for general spaces.}

In the case of Sobolev mappings between Riemannian manifolds
$W^{1,p}(M,N)$ it is
often very important to have a Riemannian metric on $N$ fixed as
changing the metric could result in turning a given variational
problem to a different one. In the case of mappings into metric
spaces $W^{1,p}(M,X)$ it is often that there is no canonical way to
define the metric in $X$ and there are several possibilities which
are in some sense equally good. One of the requirements is,
however, that the metrics be bi-Lipschitz equivalent. Actually,
from the point of view of analysis on metric spaces, spaces which
are bi-Lipschitz homeomorphic are equivalent. A fundamental
example in this context is the Heisenberg group $H^m$.
There are two canonical metrics on the Heisenberg group. The
Carnot--Carath\'eodory metric $d_C$ and the gauge metric $d$. They
are bi-Lipschitz equivalent i.e. there is a constant $C>0$ such
that for all $p,q\in H^m$
\begin{equation}
\label{equiv}
C^{-1}d_{C}(p,q) \leq d(p,q) \leq C d_{C}(p,q).
\end{equation}
Capogna and Lin, \cite{capognal} study Sobolev mappings from a
domain $\Omega$ in the Euclidean space into the Heisenberg group.
They write:
{\em because of (\ref{equiv}) both the definitions of the Lebesgue
space $L^p(\Omega,H^m)$, and the Sobolev space
$W^{1,p}(\Omega,H^m)$ do not depend on the choice between the
Carnot--Carath\'eodory distance and the one coming from the gauge
norm. To simplify the computations we will deal with the metric
$d(p,q)$.}
They consider the space of Sobolev mappings with values into the
Heisenerg group as a set and they {\em do not} equip this space
with any metric.
They are perfectly right when they say that the
set does not depend on the choice of the metric in the target as
long as the metrics that are considered are bi-Lipschitz
equivalent, like those in (\ref{equiv}).
However, if one equips the space of the Sobolev
mappings $W^{1,p}(M,X)$ from a manifold
into a metric space with a metric that comes from the
isometric embedding of $X$ into a Banach space (as was done
in \cite{heinonenkst}) it is not clear how changing the metric in
$X$ to a bi-Lipschitz equivalent one would affect the metric structure of the
space
$W^{1,p}(M,X)$.
Our paper is devoted to study this problem.

Easy examples show that different isometric embeddings of $X$
into Banach spaces may give different metrics in
$W^{1,p}(M,X)$\footnote{$\bbbr^2$ equipped with the norm
$\Vert (x_1,x_2)\Vert_1=|x_1|+|x_2|$ is a Banach space and
$i:\bbbr\to(\bbbr^2,\Vert\cdot\Vert_1)$, $i(x)=(x,0)$ if $x\geq 0$,
$i(x)=(0,x)$ if $x<0$ is an isometric embedding. It is easy to see
that this embedding equips the space $W^{1,p}([0,1],\bbbr)$ with a
metric which is different than the standard one.}
and it is not even clear what is the answer to the following
question.

\medskip

\noindent
{\bf Question 1.} {\em Does the answer to the question about the
density of Lipschitz mappings in $W^{1,p}(M,X)$ depend on the isometric
embedding of $X$ into a Banach space?}

\medskip

It is easy to see that in the case of isometric embeddings of $X$
into Hilbert spaces the answer to Question~1 is no. Since we will
be concerned with metric spaces $X$ that admit isometric
embeddings to Euclidean spaces, we will be concerned with
Question~2 (formulated below) rather than with Question~1 which is
obvious in this case. As we assume $X$ to have an isometric
embedding to the Euclidean space, we can actually identify it with
a subset of $\bbbr^\nu$ and define
\begin{equation}
\label{719}
W^{1,p}(M,X) =
\{
u\in W^{1,p}(M,\bbbr^\nu):\, \text{$u(x)\in X$ a.e.}
\}.
\end{equation}
The space $W^{1,p}(M,X)$ is equipped with the metric of the Sobolev norm
$W^{1,p}$.

Two closed manifolds can be identified if they are diffeomorphic. In the case
of metric spaces a natural counterpart of a diffeomorphism is a bi-Lipschitz
homeomorphism. This is the class of the ``most smooth'' homeomorphisms that can
be considered for arbitrary metric spaces.

Many of the metric spaces $X$ admit a bi-Lipschitz embedding into the
Euclidean space $\bbbr^\nu$, rather than an isometric one.
In such a case one could try to use this embedding
to define the space $W^{1,p}(M,X)$ by (\ref{719}). How does the space
$W^{1,p}(M,X)$ depend on the bi-Lipschitz embedding of $X$?

\medskip

\noindent
{\bf Question 2.} {\em Assume that $X$ and $Y$ are compact subsets
of $\bbbr^\nu$ that are bi-Lipschitz homeomorphic. Assume that $M$ is a
smooth closed $n$-dimensional manifold
and that Lipschitz mappings $\lip(M,X)$ are dense in
$W^{1,p}(M,X)$ for some $1\leq p<\infty$.
Are the Lipschitz mappings $\lip(M,Y)$ dense in $W^{1,p}(M,Y)$?}

\medskip

We will examine the case of Sobolev mappings
$W^{1,n}(M,X)$ from an $n$-dimensional  closed manifold
to a compact set $X\subset\bbbr^\nu$.
We will prove that if $X$ is a Lipschitz neighborhood retract
then
Lipschitz mappings $\lip(M,X)$ are dense in $W^{1,n}(M,X)$, see
Theorem~\ref{dens}.
Now if $Y$ is bi-Lipschtz homeomorphic to $X$, then $\lip(M,Y)$ is
dense in $W^{1,n}(M,Y)$ because the class of Lipschitz neighborhood
retracts is closed under bi-Lipschtz homeomorphisms
(Lemma~\ref{LNR}) and hence $Y$ is a Lipschitz neighborhood retract
as well. On the other hand
we will show in Theorem~\ref{main}
that there are two bi-Lipschitz homeomorphic compact
subsets $X$ and $Y$ of $\bbbr^{n+2}$ such that
for any $n$-dimensional closed manifold $M$
Lipschitz mappings $\lip(M,X)$ are dense in $W^{1,n}(M,X)$, but
Lipschitz mappings $\lip(M,Y)$ are not dense in $W^{1,n}(M,Y)$.
Obviously $X$ cannot be a Lipschitz neighborhood retract.

Other new results are Theorems~\ref{one},~\ref{subaru} and~\ref{smooth}.
In all theorems but Theorem~\ref{one} we deal with the Sobolev space
$W^{1,n}(M,X)$, where $n=\dim M$. One can ask similar questions for spaces
$W^{1,p}(M,X)$ with $p\neq n$.
I believe that our counterexamples
have counterparts for other exponents as well. It would be, however, more
interesting if one could prove positive results for some $p\neq n$.

\medskip

\noindent
{\bf Statement of main results.}
Consider the class of Sobolev mappings $W^{1,p}(\Om,\bbbr^\nu)$ and a bounded
Lipschitz function
$\vi\in \lip(\bbbr^\nu)\cap L^{\infty}(\bbbr^\nu)$.
Here $\Om\subset\bbbr^n$ is an open set with finite measure, but
we could also consider a closed manifold $M$ instead of $\Om$.
It is well known and easy to prove that the nonlinear operator
\begin{equation}
\label{tuptus}
\Phi:W^{1,p}(\Om,\bbbr^\nu) \to \ W^{1,p}(\Om),
\qquad
\Phi(u)=\vi\circ u
\end{equation}
is bounded i.e. it maps bounded sets to bounded sets. More precisely the
inequality $\Vert\Phi(u)\Vert_{1,p} \leq C(1+\Vert u\Vert_{1,p})$
is satisfied. The same
holds in the case in which $\Om$ is replaced by $M$.

\medskip

\noindent
{\bf Question 3.} {\em Is the operator $\Phi$ continuous?}

\medskip

This is a crucial problem for us.
Indeed, continuity  would immediately imply that $\Phi:W^{1,p}(M,X)\to
W^{1,p}(M,Y)$ is a homeomorphism whenever $\vi:X\to Y$ is a bi-Lipschitz
homeomorphism and $\Phi(u)=\vi\circ u$.
Marcus and Mizel \cite{marcusm} proved continuity of $\Phi$
given by (\ref{tuptus}) in
the case $\nu=1$. We are, however, interested mainly in the case
in which $\nu>1$.
Let us start with a positive result that will be needed later.
\begin{lemma}
\label{jeden}
Let $M$ be a smooth closed manifold and
$\vi$ a bounded Lipschitz function on $\bbbr^\nu$ that is of
the class $C^1$ away from one point $x_0$. Then the operator
$\Phi:W^{1,p}(M,\bbbr^\nu)\to W^{1,p}(M)$
is continuous for all $1\leq  p<\infty$.
\end{lemma}

Unfortunately, in general, the answer to Question~3 is in
the negative when $\nu>1$. This is shown by the following result.
I believe the result is known, but I
was unable to locate it in the literature.
In any case the proof is very easy and the phenomenon
is so important for our paper that we provide a short proof.
\begin{theorem}
\label{one} There is a Lipschitz function $\vi\in\lip(\bbbr^2)$
with compact support such that the bounded operator
$\Phi:W^{1,p}([0,1],\bbbr^2)\to W^{1,p}([0,1])$ defined as
composition $\Phi(u)=\vi\circ u$ is not continuous for any $1\leq
p<\infty$.
\end{theorem}
Similar claim follows from Theorem~\ref{main} below. Indeed, the mapping
$\Phi:W^{1,n}(M,X)\to W^{1,n}(M,Y)$ induced by the bi-Lipschitz homeomorphism
$\vi:X\to Y$ from Theorem~\ref{main} cannot be continuous. However, the proof
of Theorem~\ref{one}  is much simpler than that of Theorem~\ref{main}.

Even if we know that $\Phi:W^{1,p}(M,X)\to W^{1,p}(M,Y)$ is not
continuous for a bi-Lipschitz homeomorphism $\vi:X\to Y$
it does not immediately imply that the
spaces $W^{1,p}(M,X)$ and $W^{1,p}(M,Y)$ must have different properties.
In fact it easily follows that if $\vi:X\to Y$
is a bi-Lipschitz homeomorphism then the corresponding mapping $\Phi$ gives a
one-to-one correspondence between the spaces
$W^{1,p}(M,X)$ and $W^{1,p}(M,Y)$ and one-to-one correspondence between
the spaces $\lip(M,X)$ and $\lip(M,Y)$. Hence one
could expect positive answer to Question~2.

The following result answers Question~2 in the positive in the case in which
$p=\dim M$ and $X$ is a Lipschitz neighborhood retract.
\begin{theorem}
\label{dens}
Let $M$ be a smooth closed $n$-dimensional manifold and
$X\subset\bbbr^\nu$
a compact Lipschitz neighborhood retract. Then Lipschitz mappings
$\lip(M,X)$ are dense in $W^{1,n}(M,X)$.
\end{theorem}
Recall that $X\subset\bbbr^\nu$ is a
{\em Lipschitz neighborhood retract} if there is an open set
$X\subset U\subset\bbbr^\nu$ and a Lipschitz mapping (called {\em
Lipschitz retraction}) $p:U\to X$ such that $p(x)=x$ for all
$x\in X$. The following result is well known.
\begin{lemma}
\label{LNR}
If $X\subset\bbbr^\nu$ is a Lipschitz neighborhood retract and
$Y\subset \bbbr^\nu$ is bi-Lipschiz homeomorphic to $X$, then $Y$
is a Lipschitz neighborhood retract as well.
\end{lemma}
{\em Proof.}
If $p:U\to X$ is a Lipschitz retraction, $f:Y\to X$ a bi-Lipschitz
homeomorphism and $F:\bbbr^\nu\to\bbbr^\nu$ is a Lipschitz
extension (McShane) of $f$, then
$f^{-1}\circ p\circ F: F^{-1}(U)\to Y$ is a Lipschitz retraction.
\hfill $\Box$

\medskip

We say that $X\subset\bbbr^\nu$ is an $m$-dimensional
{\em Lipschitz submanifold} if every point $x\in X$ has a
neighborhood in $X$ which is bi-Lipschtz homeomorphic to
$B^m(0,1)$.

Theorem~\ref{dens} applies to Lipschitz submanifolds because of
the following result.
\begin{lemma}(\cite[Theorem~5.13]{LV})
Lipschitz submanifolds of $\bbbr^\nu$ are Lipschitz neighborhood
retracts.
\end{lemma}
{\em Proof.} We sketch the proof which is based on a localization
of the argument of Almgren \cite{almgren}
(cf.\ \cite[Proposition~2.13]{heinonen}).
Let $X\subset\bbbr^\nu$
be a Lipschitz submanifold. We decompose $\bbbr^\nu\setminus X$
into Whitney cubes and define the Lipschitz retraction first on
the vertices of the cubes as the projection to a closest point in
$X$. We take only those Whitney cubes which are in a close
neighborhood of $X$, as we require that all the vertices of each
cube are mapped into a domain in $X$ which is bi-Lischitz
homeomorphic to a ball. Next we extend the Lipschitz projection to
higher dimensional edges of the cubes by induction.
\hfill $\Box$

\medskip

Schoen and Uhlenbeck \cite{schoenu1}, \cite{schoenu2} proved
Theorem~\ref{dens} in the case in which $X$ is a compact smooth
manifold. They proved that the convolution approximation of $u\in
W^{1,n}(M,X)$ takes values into a tubular neighborhood of $X$ in $\bbbr^\nu$
and hence one can obtain desired approximation by composing
the convolution approximation with the smooth nearest point
projection. This argument cannot be directly applied in the case
of Theorem~\ref{dens}, even when $X$ is a Lipschitz manifold,
because now the retraction of a neighborhood of $X$ onto $X$ is
Lipschitz only and according to Theorem~\ref{one} the composition
with a Lipschitz function need not be continuous in the Sobolev
norm.
Proof of Theorem~\ref{dens} will require a much more delicate
approximation method than the one coming from convolution.

A nice and nontrivial class of examples of Lipschitz manifolds was
provided by Toro, \cite{toro} (see also \cite{mullers}).
She proved that surfaces in $\bbbr^3$ that are locally graphs of
$W^{2,2}$ functions are Lipschitz submanifolds.
If we replace the class of $W^{2,2}$ surfaces by
continuous $W^{1,2}$ surfaces, then we can construct a
counterexample to the density. Such a counterexample is provided
by the case $n=2$ of the following result.

\begin{theorem}
\label{subaru}
There is a subset $N\subset\bbbr^{n+1}$ homeomorphic to $S^n$ and
$x_0\in N$ such that
$N\setminus\{ x_0\}$ is a smooth submanifold of $\bbbr^{n+1}$
and $N$ is a graph of
a continuous $W^{1,n}$ function near $x_0$. That means $N$ is a
regular $n$-dimensional submanifold of $\bbbr^{n+1}$ with one
singular point. This construction can be done in such a way that
for every closed $n$-dimensional manifold $M$, Lipschitz mappings
$\lip(M,N)$ are not dense in $W^{1,n}(M,N)$.
\end{theorem}

One should compare Theorem~\ref{subaru} with Theorem~\ref{dens}.
As we will see Theorem~\ref{subaru} is a consequence of the proof
of Theorem~\ref{smooth}. Namely we will see that if
$M=S^n\subset\bbbr^{n+1}$, then the singular manifold
$N=\widetilde{M}$ constructed in the proof of Theorem~\ref{smooth}
has properties required in Theorem~\ref{subaru}.

If $X\subset\bbbr^\nu$ is a subset,
then by smooth mappings $C^{\infty}(M,X)$ we will mean
the class of all mappings $u\in C^{\infty}(M,\bbbr^\nu)$ such that
$u(M)\subset X$.
\begin{theorem}
\label{main}
Fix an integer $n\geq 2$. There is a compact and connected
set $X\subset\bbbr^{n+2}$ and a global
bi-Lipschitz homeomorphism
$\Phi:\bbbr^{n+2}\to\bbbr^{n+2}$ with the property
that for any closed $n$-dimensional manifold $M$
smooth mappings $C^{\infty}(M,X)$
are dense in $W^{1,n}(M,X)$, but Lipschitz mappings $\lip(M,Y)$ are not dense
in $W^{1,n}(M,Y)$, where $Y=\Phi(X)$.
\end{theorem}
According to Theorem~\ref{dens} and Lemma~\ref{LNR}, $X$ cannot be
a Lipschitz neighborhood retract. In particular
there is no such example as in Theorem~\ref{main} in the case in
which $X$ is a closed manifold. However, if we allow a global homeomorphism
$\Phi$ of
the Euclidean space to be locally bi-Lipschitz everywhere but in one point,
then one can find counterexamples even in the
class of closed manifolds.
\begin{theorem}
\label{smooth}
Let $M\subset\bbbr^{\nu}$ be a closed $n$-dimensional manifold.
Then there is a homeomorphism
$\Phi\in C^{\infty}(\bbbr^\nu,\bbbr^\nu)$ which is a
diffeomorphism in $\bbbr^{\nu}\setminus\{ 0\}$,
which is identity outside a
sufficiently large ball and has the property that
Lipschitz mappings $\lip(M,\widetilde{M})$ are not dense in
$W^{1,n}(M,\widetilde{M})$, where $\widetilde{M}=\Phi^{-1}(M)$.
\end{theorem}
The example constructed in the proof of Theorem~\ref{smooth}
will be employed in the proof of Theorem~\ref{main}.

\medskip

\noindent
{\bf Notation.}
Notation employed in the paper is rather standard.
The $k$-dimensional Hausdorff measure will be denoted by $\H^k$
and the average value by
$$
u_E = \mvint_{E} u\, d\mu = \mu(E)^{-1} \int_{E} u\, d\mu\, .
$$
Symbol $B$ will be used to denote a ball. By $\sigma B$,
where $\sigma\geq 1$ we will denote a ball concentric with $B$ and with
the radius $\sigma$ times that of $B$.
We say that $M$ is a closed manifold if it is smooth, compact and
without boundary.
The $L^p$ norm of a function $u$ will be denoted by $\Vert u\Vert_{p}$.
The Sobolev space $W^{1,p}(M)$ consists of all $u\in L^{p}(M)$ such that
$\nabla u\in L^{p}(M)$. The space is equipped with the norm
$\Vert u\Vert_{1,p}=\Vert u\Vert_{p} + \Vert \nabla u\Vert_{p}$.
By $C$ we will denote a general constant
which can change its value even in the same string of estimates.
We will write $f:X\twoheadrightarrow Y$ to designate
any function $f:X\to Y$ such
that $f(X)=Y$. Coordinates in $\bbbr^n$ will often be denoted by $(x',x^n)$,
where $x'\in \bbbr^{n-1}$.

\medskip

\noindent
{\bf Acknowledgement.}
I wish to express my deepest gratitude to the referee
for the incredibly careful work. Referee's comments helped me a
great deal to improve Theorem~\ref{dens} and other parts of the
paper as well.

\section{Proof of Lemma~\ref{jeden}}

For $f_n\to f$ in $W^{1,p}(M,\bbbr^\nu)$ and a.e. we will show
that
$\Phi\circ f_n\to\Phi\circ f$ in $W^{1,p}$.
We split $M$ into three pairwise disjoint sets
$A=f^{-1}(x_0)$,
$$
B_k=\left( \bigcup_{n=k}^{\infty} f_{n}^{-1}(x_0)\right)
\setminus f^{-1}(x_0)
\qquad
C_k= M\setminus
\left( f^{-1}(x_0) \cup \bigcup_{n=k}^{\infty}
f_{n}^{-1}(x_0)\right).
$$
Note that the measure of $B_k$ goes to $0$ as $k\to\infty$ and
hence for every $\eps>0$ there is $k$ such that
$$
\int_{B_k} |D(\Phi\circ f)|^p +
\sup_{n\geq k} \int_{B_k} |D(\Phi\circ f_n)|^p <\eps,
$$
because the family
$|D(\Phi\circ f_n)|^p\leq C|D f_n|^p$ is equiintegrable as
$D f_n\to Df$ in $L^p$.
In particular
$$
\int_{B_k} |D(\Phi\circ f_n)-D(\Phi\circ f)|^p < C\eps
$$
for all $n\geq k$.
The chain rule applies to
$D(\Phi\circ f_n)$ and to $D(\Phi\circ f)$ on $C_k$ when $n\geq k$
and hence $D(\Phi\circ f_n)\to D(\Phi\circ f)$ in $L^p$ on $C_k$.
Sine $f$ is constant on $A$, $Df=0$ a.e. on $A$ and $D(\Phi\circ
f)=0$ a.e. on $A$. Hence $Df_n\to 0$ in $L^p$ on $A$. Thus
$|D(\Phi\circ f_n)|\leq C|Df_n|\to 0$ in $L^p$ on $A$. This,
however, implies that $D(\Phi\circ f_n)\to 0=D(\Phi\circ f)$ in
$L^p$ on $A$. All the facts put together easily imply that
$D(\Phi\circ f_n)\to D(\Phi\circ f)$ in $L^p(M)$. The proof is
complete.
\hfill $\Box$

\section{Proof of Theorem~\ref{one}}

We will construct  a bounded Lipschitz function $\vi:\bbbr^2\to\bbbr$ with
compact support $\supp \vi =[0,2]\times [0,1]$ such that the bounded operator
\[
\Phi:W^{1,p}([0,1],\bbbr^2) \to W^{1,p}([0,1],\bbbr),
\qquad
\Phi(u)=\vi\circ u
\]
is not continuous for any choice of $p$ in the interval $1\leq p<\infty$.

To this end we first define a compact set $K\subset [0,2]\times [0,1]$ as
union of an infinite family of segments defined below.

Let $a_i=2-2^{-i+1}$ for $i=0,1,2,\ldots$ We define
\[
I_i=\{ a_i\} \times [0,1]
\quad
\mbox{and}
\quad
I_\infty =\{ 2\} \times [0,1].
\]
This is a family of vertical segments. Next for each $i=0,1,2,\ldots$ we define
a finite family of horizontal segments $\{ J_{i,k}\}_{k=0}^{2^i}$
that join $I_i$ with $I_{i+1}$ Namely we set
\[
J_{i,k} = [a_i,a_{i+1}]\times \{ k/2^i\}.
\]
Now we define the set $K$ as follows
\[
K= I_\infty \cup
\left( \bigcup_{i=0}^{\infty} I_i\right) \cup
\left( \bigcup_{i=0}^{\infty} \bigcup_{k=0}^{2^i} J_{i,k} \right).
\]
The set $K$ is the union of boundaries an infinite family of squares. Now we
define the function $\vi$ by the formula
$$
\vi(x) =
\left\{
\begin{array}{ccc}
\dist(x,K)    & \mbox{if $x\in [0,2]\times [0,1]$,}\\
0 &   \mbox{if $x\in\bbbr^2\setminus [0,2]\times [0,1]$.}
\end{array}
\right.
$$
The graph of $\vi$ over each of the squares has the shape of a
pyramid, so the graph of $\vi$ is an infinite collection of
pyramids.
The Lipschitz constant of $\vi$ is $1$. Fix $1\leq p<\infty$. We will show that
the operator $\Phi$ is not continuous in $W^{1,p}$. Let  $u_i:[0,1]\to\bbbr^2$,
$u_i(x) =((a_i+a_{i+1})/2,x)$. The image of $u_i$ is the vertical segment
that is in the middle between $I_i$ and $I_{i+1}$. Obviously $u_i\to u$ in
$W^{1,p}$ as $i\to\infty$, where $u:[0,1]\to\bbbr^2$,
$u(x)=(2,x)$. However, $\vi\circ u_i$ does not
converge to $\vi\circ u\equiv 0$ in $W^{1,p}$ because
$(\vi\circ u_i)' =\pm 1$ a.e.
The proof is complete.
\hfill $\Box$

\section{Proof of Theorem~\ref{dens}.}

The approximation technique
employed here is similar to that used in \cite{bojarskihs}, \cite{hajlaszki}.
Fix a Riemannian tensor on the manifold $M$. Denote the induced
distance
and the gradient by $d(x,y)$ and $\nabla u$ respectively.
$\M g(x) = \sup_{r>0} \mvint_{B(x,r)}|g(z)|\, dz$ will be used to denote
the Hardy-Littlewood maximal function.

There exists $r_0$ depending on the Riemannian structure of $M$
such that the following inequalities
$$
\mvint_{B}|u-u_B| \leq C r \mvint_{B} |\nabla u|,
$$
and
\begin{equation}
\label{x361}
|u(x) - u_B| \leq C r \M |\nabla u|(x)
\qquad
\text{for a.e. $x\in B$}.
\end{equation}
hold true whenever the radius $r$ of $B$ is less than $r_0$ and
$u\in W^{1,p}(B)$, $1\leq p<\infty$. The first inequality is known
as the Poincar\'e inequality.
The second inequality follows e.g. from Lemma~1.50 and
Theorem~1.32(i) with $\alpha=1$ in \cite{malyz}.

As a corollary one can easily prove the following version of a
well known pointwise inequality, see e.g.
\cite{acerbif}, \cite{hajlasz3},
\cite[Theorems~3.2~and~3.3]{hajlaszk} and references therein.
\begin{lemma}
\label{pointwise}
Let $M$ be a closed manifold and $u\in W^{1,p}(M)$, where $1\leq p<\infty$.
Then
$$
|u(x)-u(y)| \leq C d(x,y)(\M |\nabla u|(x) + \M |\nabla u|(y))
\qquad
\text{for a.e. $x,y\in M$}.
$$
\end{lemma}

Let $u\in W^{1,n}(M,X)$. Our aim is to construct a family of
Lipschitz mappings $u_t\in\lip(M,\bbbr^\nu)$ such that
\begin{enumerate}
\item[(A)] The Lipschitz constant of $u_t$ is bounded by $Ct$.
\item[(B)] $t^n \H^n(\{ u\neq u_t\})\to 0$ as $t\to\infty$.
\item[(C)] $\sup_{x\in M} \dist(u_t(x),X)\to 0$ as $t\to\infty$.
\end{enumerate}
Before we construct a family $\{u_t\}$ we will show how to use
the properties (A), (B) and (C) to complete the proof of the
theorem.

Let $p:U\to X$ be a Lipschitz neighborhood retraction. By (C),
$u_t(M)\subset U$ for all sufficiently large $t$ and hence
$p\circ u_t:M\to X$ is a family of Lipschitz mappings well defined
for large $t$. We will show now that $p\circ u_t\to u$ in
$W^{1,n}$ when $t\to\infty$. Indeed,
$p\circ u_t\to u$ a.e. by (B) and hence $p\circ u_t\to u$ in
$L^n$ (mappings are bounded). Finally (A) and (B) imply
\begin{eqnarray*}
\lefteqn{
\int_{M} |D(p\circ u_t)-Du|^{n}
=
\int_{ \{p\circ u_t\neq u\} } |D(p\circ u_t)-Du|^{n} } \\
& \leq &
C t^n \H^{n}(\{p\circ u_t\neq u\}) + C \int_{ \{ p\circ u_t\neq u\} } |Du|^{n}
\to 0,
\qquad
\mbox{as $t\to\infty$}
\end{eqnarray*}

Now we proceed to the construction of the family $\{u_t\}$. Let
$E_t=\{ x\in M:\, \M |\nabla u|(x)\leq t\}$.
The set $M\setminus E_t$ is open.
There is a Whitney decomposition of $M\setminus E_t$ into balls and
subordinated Lipschitz partition of unity. This is to say there is a constant
$C\geq 1$ depending on the Riemannian structure of $M$ only and a sequence $\{
x_i\}_{i\in I}$ of points in $M\setminus E_t$ such that with
$r_i=\dist(x_i,E_t)/10$ we have
\begin{enumerate}
\item[(a)] $\bigcup_{i\in I} B(x_i,r_i)= M\setminus E_t$;
\item[(b)] $B(x_i,5r_i)\subset M\setminus E_t$ for all $i\in I$;
\item[(c)] for every $i\in I$ and all $x\in B(x_i,5r_i)$ we have
$5 r_i \leq \dist (x,E_t) \leq 15 r_i$;
\item[(d)] no point of $M\setminus E_t$ belongs to more than $C$ balls
$\{B(x_i,5r_i)\}_{i\in I}$;
\item[(e)] There is a family of Lipschitz continuous functions
$\{\vi_i\}_{i\in I}$ such that $\supp\vi_i\subset B(x_i,2r_i)$,
$0\leq \vi_i\leq 1$, $\sum_{i\in I}\vi_i= 1$ and the Lipschitz constant of
$\vi_i$ is bounded by $Cr_i^{-1}$.
\end{enumerate}

The proof of the existence of
a Whitney decomposition of a domain in $\bbbr^n$ into cubes
with pairwise disjoint interiors can be found in \cite{stein}. Standard
modification of the method gives the above decomposition into balls.
Actually the construction is  so general that it can be carried out on
arbitrary doubling metric
spaces, see \cite{coifmanw}, \cite[Lemma~2.9]{maciass}.
In what follows we will use notation
$B_i=B(x_i,r_i)$ and $2B_i=B(x_i,2r_i)$.
We define
\[
u_t =
\left\{
\begin{array}{ccc}
u(x)    & \mbox{for $x\in E_t$,}\\
\sum_{i\in I} \vi_i(x) u_{B_i} &   \mbox{for $x\in M\setminus E_t$.}
\end{array}
\right.
\]
Note that by Lemma~\ref{pointwise}
the function $u_t|_{E_t}=u|_{E_t}$ is Lipschitz continuous with the Lipschitz
constant $Ct$.

By the Hardy--Littlewood theorem, \cite{stein}, $\M |\nabla u|\in
L^n$, and hence $t^n\H^n(M\setminus E_t)\to 0$ as $t\to\infty$.
This immediately implies the property (B) of $u_t$.

Now we prove the property (A).
Given $x\in M\setminus E_t$, let $\overline{x}\in E_t$ be such that
$d(x,\overline{x})=\dist (x,E_t)$. Note that
whenever $B_i$ is such that $x\in 2B_i$, then $B_i\subset
B(\overline{x},2d(x,\overline{x})):=\widetilde{B}(\overline{x})$. Observe that
radius of $\widetilde{B}(\overline{x})$ is comparable to $r_i$ whenever $x\in
2B_i$.
Choose $t$ large enough to make the measure $\H^n(M\setminus E_t)$
so small that radii of all of the balls
$\widetilde{B}(\overline{x})$ are less than $r_0$.
For $x\in M\setminus E_t$ we have
\begin{equation}
|u_t(x)-u_t(\overline{x})|
 =
\left| \sum_{i\in I} \vi_{i}(x)(u(\overline{x}) -u_{B_i})\right|
\leq
\sum_{i\in I_x} |u(\overline{x})-u_{B_i}|,
\label{last}
\end{equation}
where $I_x=\{i\in I:\, x\in\supp\vi_i\}$
(observe that $x\in 2B_i$ for $i\in I_x$).
We will prove now that
\begin{equation}
\label{asia}
|u(\overline{x}) - u_{B_i}| \leq C t r_i
\qquad
\mbox{for all $i\in I_x$.}
\end{equation}
This together with inequality (\ref{last}), the observation
that $r_i\leq d(x,\overline{x})$ for $i\in I_x$ and the fact that the number of
elements in $I_x$ is bounded by a constant that does not depend on $x$
will give
\begin{equation}
\label{trzy}
|u_t(x) -u_t(\overline{x})| \leq
C t d(x,\overline{x}).
\end{equation}
To prove (\ref{asia}) note that
$|u(\overline{x}) - u_{B_i}| \leq
|u(\overline{x}) - u_{\widetilde{B}(\overline{x})}| +
|u_{\widetilde{B}(\overline{x})} -u_{B_i}|$.
Applying inequality (\ref{x361}) yields
\[
|u(\overline{x}) - u_{\widetilde{B}(\overline{x})}|
\leq
C r_i \M |\nabla u|(\overline{x}) \leq C t r_i
\]
because $\overline{x}\in E_t=\{\M |\nabla u|\leq t\}$.
To estimate the second term we apply
the Poincar\'e inequality
\begin{eqnarray*}
|u_{\widetilde{B}(\overline{x})} -u_{B_i}|
& \leq &
\mvint_{B_i} |u-u_{\widetilde{B}(\overline{x})}|
\leq
C\mvint_{\widetilde{B}(\overline{x})} |u-u_{\widetilde{B}(\overline{x})}| \\
& \leq &
C r_i \mvint_{\widetilde{B}(\overline{x})} |\nabla u|
\leq
C r_i \M |\nabla u|(\overline{x})
\leq
C t r_i.
\end{eqnarray*}
The above inequalities complete the proof of (\ref{asia}) and hence that for
(\ref{trzy}). We proved inequality (\ref{trzy}) for $x\in M\setminus E_t$. The
inequality holds, however, for all $x\in M$ because $x=\overline{x}$
whenever $x\in E_t$.

Now we are ready to estimate the Lipschitz constant of $u_t$.
If $d(x,y)\geq \min\{ \dist(x,E_t), \dist(y,E_t)\}$, then
\begin{eqnarray*}
|u_t(x) - u_t(y)|
& \leq &
|u_t(x) - u_t(\overline{x})| +
|u_t(\overline{x}) - u_t(\overline{y})| +
|u_t(\overline{y}) - u_t(y)| \\
& \leq &
C t (d(x,\overline{x}) + d(\overline{x},\overline{y}) + d(\overline{y},y))
\leq
C t d(x,y).
\end{eqnarray*}
We employed here the fact that $u_t|_{E_t} = u|_{E_t}$ is Lipschitz continuous
with the Lipschitz constant $Ct$.
Thus we can assume that
$d(x,y) \leq \min\{\dist(x,E_t),\dist(y,E_t)\}$.
Since $\sum_{i\in I}(\vi_i(x)-\vi_i(y))=0$ we have
\begin{eqnarray*}
\lefteqn{
|u_t(x) - u_t(y)|
 =
\left| \sum_{i\in I} \vi_i(x) u_{B_i} -
\sum_{i\in I} \vi_i(y) u_{B_i} \right|} \\
& = &
\left| \sum_{i\in I} (\vi_i(x) - \vi_i(y))(u(\overline{x})-u_{B_i})\right|
 \leq
C d(x,y) \sum_{i\in I_x\cup I_y} r_i^{-1}|u(\overline{x})-u_{B_i}| \\
& \leq &
Cd(x,y) \sum_{i\in I_x\cup I_y} r_i^{-1} Ctr_i
\leq C t d(x,y).
\end{eqnarray*}
In the inequality second to last we employed the fact that
$|u(\overline{x})-u_{B_i}|\leq Ctr_i$ for all $i\in
I_x\cup I_y$. The proof of this fact is very similar to the proof of inequality
(\ref{asia}). The proof of the second claim is complete.

Now we will prove the property (C) of $u_t$.
For $x\in E_t$ we have $u_t(x) = u(x)\in X$ and hence we can assume that $x\in
M\setminus E_t$.
Given $\eps>0$ there is $r_0$ such that for every ball $B$ in $M$ with radius
$r\leq r_0$ we have
\begin{equation}
\label{ast}
\left( \int_{B} |\nabla u|^{n} \right)^{1/n} \leq \eps.
\end{equation}
Hence if we choose $t$ large enough we can guarantee that (\ref{ast}) holds for
each of the balls $B=\widetilde{B}(\overline{x})$.
Since
$\dist(u_{\widetilde{B}(\overline{x})},X) \leq
|u_{\widetilde{B}(\overline{x})}-u(z)|$ for
$z\in \widetilde{B}(\overline{x})$, the Poincar\'e
inequality yields
\begin{eqnarray}
\lefteqn{
\dist(u_{\widetilde{B}(\overline{x})},X) \leq
\mvint_{\widetilde{B}(\overline{x})} |u-u_{\widetilde{B}(\overline{x})}|}
\nonumber \\
\label{b52}
& \leq &
C r_i \mvint_{\widetilde{B}(\overline{x})} |\nabla u| \leq
C \left(\int_{\widetilde{B}(\overline{x})}|\nabla u|^{n} \right)^{1/n} \leq
C\eps.
\end{eqnarray}
This means for sufficiently large $t$ all averages
$\widetilde{B}(\overline{x})$ are very close to $X$. Since
\[
\dist(u_t(x),X) \leq |u_t(x) - u_{\widetilde{B}(\overline{x})}| + C\eps
\]
it suffices to prove that
$|u_t(x) - u_{\widetilde{B}(\overline{x})}|<C\eps$. We have
\begin{eqnarray*}
\lefteqn{
|u_t(x) - u_{\widetilde{B}(\overline{x})}|
=
\left| \sum_{i\in I} \vi_{i}(x)
(u_{B_i}- u_{\widetilde{B}(\overline{x})}) \right|} \\
& \leq &
\sum_{i\in I_x} \mvint_{B_i} |u-u_{\widetilde{B}(\overline{x})}| \leq
C \sum_{i\in I_x} \mvint_{\widetilde{B}(\overline{x})}
|u-u_{\widetilde{B}(\overline{x})}| \leq C\eps.
\end{eqnarray*}
The last inequality follows from the Poincar\'e inequality applied
as in (\ref{b52}) and the fact that the number of indices in $I_x$ is bounded
by a constant independent of $x$. This completes the proof of the last
claim and hence that of the theorem.
\hfill $\Box$

\section{Proof of Theorem~\ref{smooth}.}
\label{four}

Let us start with a brief description of the main steps in the proof.
Let $M\subset\bbbr^\nu$ be a closed $n$-dimensional manifold. In the first step
we define a continuous $W^{1,n}$ function $\gamma$ on $\bbbr^n$ with compact
support contained in the unit ball $B^n(0,1)$.
In the second step we use the fact that $M$
contains a subset diffeomorphic to $\overline{B^{n}}(0,1)$
to replace a part
of $M$ by the graph of $\gamma$. The resulting space is denoted by
$\widetilde{M}$. It is still a subset of $\bbbr^\nu$.
The function $\gamma$ is constructed in such a way that it is
smooth everywhere but in one point
and hence $\widetilde{M}$ is a smooth manifold with one nonsmooth point on it.
Obviously $\widetilde{M}$ is homeomorphic to $M$, actually
a careful construction of $\widetilde{M}$ can make it $C^\infty$
homeomorphic in the sense that there is a homeomorphism $\Phi\in
C^{\infty}(\bbbr^\nu,\bbbr^\nu)$ such that $\Phi(\widetilde{M})=M$. Moreover
$\Phi$ can be chosen to be a diffeomorphism everywhere but in one point and
identity outside a sufficiently large ball.
While $\Phi|_{\widetilde{M}}:\widetilde{M}\to M$ is a smooth
homeomorphism, it is shown that there
is no smooth homeomorphism from $M$ onto $\widetilde{M}$.
Actually due to high oscillations of
$\gamma$ near the nonsmooth point, it is proved that there is no Lipschitz
mapping $\psi:M\twoheadrightarrow\widetilde{M}$.

In the last step we show that Lipschitz mappings $\lip(M,\widetilde{M})$ are
not dense in $W^{1,n}(M,\widetilde{M})$. The reason being the following:
$\gamma\in W^{1,n}$ and hence there is a homeomorphism
$\vi:M\twoheadrightarrow\widetilde{M}$ of the class $W^{1,n}$. It follows
from the homotopy properties of $W^{1,n}$ mappings that
every Lipschitz mapping
$\psi:M\to\widetilde{M}$ sufficiently close to $\vi$ in the $W^{1,n}$ norm is
homotopic to $\vi$ and hence
$\psi(M)=\widetilde{M}$.  This contradicts, however, the lack of
Lipschitz mappings $\psi:M\twoheadrightarrow\widetilde{M}$.

Actually the steps of the proof will be organized in a slightly different
order due to the idea of collecting technical constructions at the end of the
proof.

\medskip
\noindent
{\bf Construction of $\gamma$.} The construction of
$\gamma$ presented below is closely related to the fact that the
$n$-capacity of a point is zero. Those who know what it means can
avoid the part of the construction involving truncation and convolution, but
just invoke basic definitions.

It is well known and easy to prove
that
$\eta(x)=\log|\log|x||\in W^{1,n}(B^{n}(0,e^{-1}))$.
Define the truncation of $\eta$ between levels $s$ and
$t$, $0<s<t<\infty$ by
\[
\eta_{s}^{t}(x) =
\left\{
\begin{array}{ccc}
t-s    & \mbox{if $\eta(x)\geq t$,}\\
\eta(x)-s &   \mbox{if $s\leq \eta(x) \leq t$}\\
0 & \mbox{if $\eta(x) \leq s$.}
\end{array}
\right.
\]
Fix an arbitrary $\tau>0$. For every $\eps>0$ there is a sufficiently large $s$
such that
$\widetilde{\gamma}_{\eps,\tau}:=\eta_{s}^{s+\tau}$ is a Lipschitz function on
$\bbbr^n$ with the properties:
\[
\supp\widetilde{\gamma}_{\eps,\tau}\subset B^{n}(0,\eps/2),
\]
\[
\text{$0\leq\widetilde{\gamma}_{\eps,\tau}\leq\tau$
and $\widetilde{\gamma}_{\eps,\tau}=\tau$ in a neighborhood of $0$},
\]
\[
\int_{\bbbr^n}|\nabla \widetilde{\gamma}_{\eps,\tau}|^{n} < \eps^n.
\]
The function $\widetilde{\gamma}_{\eps,\tau}$ is not smooth because it is
defined as a truncation, however, mollifying
$\widetilde{\gamma}_{\eps,\tau}$
gives a smooth function, denoted by
$\gamma_{\eps,\tau}$, with the same properties
as those of $\widetilde{\gamma}_{\eps,\tau}$ listed above.

Let $(x',x^n)$ denotes the coordinates in $\bbbr^n$,
where $x'=(x^1,\ldots,x^{n-1})$. Let
$x_k=(x_k',x_k^n)=(0,2^{-k})$, $k=6,7,8\ldots$, and let
$B_k=B^n(x_k,2^{-(k+2)})$. The balls $B_k$ are pairwise disjoint
and all contained in $B^n(0,2^{-5})$.

Observe that we can find $N_k=c(n) 2^{(k^2-k)n}$ pairwise disjoint balls of
radius $2^{-k^2}$ inside $B_k$. Select such a family for each
$B_k$ and denote the balls in the family by $B_{k,i}$, $i=1,2,\ldots,
N_k$. Denote the centers of $B_{k,i}$ by $x_{k,i}$ and define
\begin{equation}
\label{myszka}
\gamma(x) = \sum_{k=6}^{\infty} \sum_{i=1}^{N_k} \eta_{k,i}(x),
\end{equation}
where $\eta_{k,i}(x) = \gamma_{2^{-k^2}, 2^{-k}}(x-x_{k,i})$.
Each of the functions $\eta_{k,i}$
is supported in $\frac{1}{2}B_{k,i}$ and hence the
supports of  the functions in the double sum (\ref{myszka}) are pairwise
disjoint. This immediately implies that $\gamma\in
C^{\infty}(\bbbr^n\setminus\{ 0\})$, $\gamma$ is continuous on $\bbbr^n$,
$\gamma(0)=0$ and $\supp \gamma \subset B^{n}(0,2^{-5})$. Moreover
\[
\int_{\bbbr^n} |\nabla \gamma|^{n} =
\sum_{k=6}^{\infty} \sum_{i=1}^{N_k}
\int_{\bbbr^n}|\nabla \eta_{k,i}|^{n} <
\sum_{k=6}^{\infty} \sum_{i=1}^{N_k} 2^{-k^2 n} =
c(n) \sum_{k=6}^{\infty} 2^{-kn} <\infty,
\]
and hence $\gamma\in W^{1,n}(\bbbr^n)$. Denote the graph of $\gamma$ over
$\overline{B^n}(0,1)$ by $\Gamma$ i.e.
\[
\Gamma = \{(x,\gamma(x))\in\bbbr^{n+1}:\, x\in \overline{B^n}(0,1) \}.
\]
We equip $\Gamma$ with the metric induced from $\bbbr^{n+1}$. The following
lemma is a consequence of  high oscillations of $\gamma$ near the origin.
\begin{lemma}
\label{niema}
Let $N$ be a closed $n$-dimensional Riemannian manifold. Then there is no
Lipschitz mapping $\psi:N\twoheadrightarrow\Gamma$.
\end{lemma}
{\em Remark.}
The claim of the lemma remains the same (with the same proof)
if we replace $N$ by any metric space $X$ equipped with a finite Borel measure
$\mu$ satisfying the lower bound for the growth
$\mu(B(x,r))\geq C r^n$ for all $x\in X$ and all $0<r\leq r_0$. Clearly $N$
equipped with the Hausdorff measure $\H^n$ is an example.

{\em Proof of Lemma~\ref{niema}.}
We argue by contradiction. Suppose $\psi:N\twoheadrightarrow\Gamma$ is an
$L$-Lipschitz mapping. Fix a positive integer $k$. The graph of $\eta_{k,i}$
over $\frac{1}{2}B_{k,i}$ is of height $2^{-k}$ and hence we can find
$2^{k^{2}-k}$ points on that graph whose mutual distances are at least
$2^{-k^2}$. Choose such points and denote them by $y_{k,i,j}$,
$j=1,2,\ldots,2^{k^2-k}$. Since distances between the balls
$\{\frac{1}{2}B_{k,i}\}_{i}$ are no less than $2^{-k^2}$
we conclude that mutual distances
of the points in the family $\{ y_{k,i,j}\}_{i,j}$ are at least $2^{-k^2}$. The
family $\{ y_{k,i,j}\}_{i,j}$ consists of
$$
c(n) 2^{(k^2-k)n} 2^{k^2-k} = c(n) 2^{(k^2-k)(n+1)}
$$
points. Choose $z_{i,j}\in\psi^{-1}(y_{k,i,j})$ arbitrarily.
Since $\psi$ is $L$-Lipschitz
the distances between points $\{ z_{i,j}\}_{i,j}$ are at least
$L^{-1}2^{-k^2}$. Hence the balls
$B(z_{i,j},L^{-1}2^{-k^2}/2)$ are pairwise disjoint. This in turn gives the
lower estimate for the volume of $N$.
$$
\H^n(N) \geq
\sum_{i,j} \H^{n}(B(z_{i,j},L^{-1}2^{-k^2}/2)) \geq
C(L^{-1}2^{-k^{2}}/2)^{n} 2^{(k^2-k)(n+1)}\to \infty
$$
as $k\to\infty$. This clearly contradicts the finiteness of the volume of $N$.
The proof of the lemma is complete.
\hfill $\Box$

\medskip
\noindent
{\bf The construction of $\widetilde{M}$.}
Take a subset of $M$ that is diffeomorphic to $\overline{B^n}(0,1)$ and
replace it by the graph $\Gamma$. Denote the resulting space by
$\widetilde{M}\subset\bbbr^\nu$. The construction is quite obvious and does not
require further explanations. In any case, we will present all details
in the next subsection because it will be needed for the
construction of $\Phi$.

The following properties of $\widetilde{M}$ are obvious. The property (c)
below is a consequence of the fact that the mapping
\[
\Gamma\to \overline{B^n}(0,1),
\qquad
(x,\gamma(x))\mapsto x
\]
is Lipschitz continuous, and property (d) is a consequence of Lemma~\ref{niema}.
\begin{enumerate}
\item[(a)] There is $\widetilde{x}_{0}\in\widetilde{M}$ such that
$\widetilde{M}\setminus\{\widetilde{x}_0\}$ is a smooth submanifold
of $\bbbr^\nu$.
\item[(b)] There is a homeomorphism
$\vi:M\twoheadrightarrow\widetilde{M}$ in the
class $\vi\in W^{1,n}(M,\widetilde{M})$ with the property that
$\vi:M\setminus\{ x_0\}\twoheadrightarrow
\widetilde{M}\setminus\{ \widetilde{x}_0\}$ is a $C^\infty$ diffeomorphism of
manifolds, where $x_0=\vi^{-1}(\widetilde{x}_0)$.
\item[(c)] $\vi$ can be choosen in such a way that
$\vi^{-1}:\widetilde{M}\twoheadrightarrow M$ is Lipschitz continuous.
\item[(d)] There is no Lipschitz mapping
$\psi:M\twoheadrightarrow\widetilde{M}$.
\end{enumerate}
The above properties do not immediately
guarantee the existence of a homeomorphism $\Phi\in
C^{\infty}(\bbbr^\nu,\bbbr^\nu)$ being diffeomorphism
everywhere but in one point and such
that $\Phi(\widetilde{M})=M$. One could try to construct $\Phi$ as an
extension of $\vi^{-1}:\widetilde{M}\twoheadrightarrow M$, however,
high oscillations of $\widetilde{M}$ near $\widetilde{x}_0$ seem to produce
a singularity of $D\Phi$ at $\widetilde{x}_0$, and this is just one of the
technical problems that we will have to face.

We postpone the technically complicated
construction of $\Phi$  for a while and show that
Lipschitz mappings $\lip(M,\widetilde{M})$ are not dense in
$W^{1,n}(M,\widetilde{M})$.

We will actually show that the homeomorphism
$\vi\in W^{1,n}(M,\widetilde{M})$ cannot be approximated by
Lipschitz mappings $\lip(M,\widetilde{M})$ in the $W^{1,n}$ norm. Suppose by
contrary that $\vi_k:M\to\widetilde{M}$ is a sequence of Lipschitz mappings
such that $\vi_k\to \vi$ in $W^{1,n}$. Then
$\vi^{-1}\circ\vi_k\to\vi^{-1}\circ\vi=\id:M\to M$ in $W^{1,n}(M,M)$.
Indeed, by Lemma~\ref{jeden}
the composition with $\vi^{-1}$ is continuous in the $W^{1,n}$
norm.\footnote{
Note that in the construction that will be carried out later
$\vi^{-1}$ admits an extension to
$\Phi\in C^{\infty}(\bbbr^\nu,\bbbr^\nu)$.
Hence $\vi^{-1}$ is
smooth everywhere and there is no problem with the continuity of the operator
of composition.}.
We will also need the following result of White~\cite{white}.
\begin{lemma}
\label{5.2}
Let $M$ and $N$ be closed manifolds and let $n=\dim M$. Then there exists
$\eps>0$ such that any two Lipschitz mappings $g_1,g_2:M\to N$ satisfying
$\Vert g_1-g_2\Vert_{1,n}<\eps$ are homotopic.
\end{lemma}
The degree modulo $2$ of the
identity mapping ${\rm id}:M\to M$ equals
$1$ (we do not assume $M$ be orientable, so the
Brouwer degree needs to be replaced by degree modulo $2$, see
\cite{milnor}).
Lemma~\ref{5.2} implies that for sufficiently large $k$,
$\vi^{-1}\circ\vi_k:M\to M$ is homotopic to the identity map.
Hence by the Homotopy Lemma \cite{milnor},
the degree modulo $2$ of $\vi^{-1}\circ\vi_k$ equals $1$ as well,
and the definition of the degree modulo $2$ yields
$(\vi^{-1}\circ \vi_k)(M)=M$.
This, in turn, implies $\vi_k(M)=\vi(M)=\widetilde{M}$ which
contradicts property (d).

\medskip

\noindent
{\bf The construction of $\Phi$.}
To complete the proof of the theorem it remains to show that there is a
$C^\infty$ homeomorphism $\Phi:\bbbr^\nu\to\bbbr^\nu$ which is diffeomorphism
everywhere but in one point, identity outside a sufficiently large
ball, and satisfies $\Phi(\widetilde{M})=M$.
To this end we will repeat the construction of
$\widetilde{M}$ providing more details.
In what follows $B^k(0,1)$, $1\leq k<\nu$, will be identified with
a subset of $\bbbr^\nu$ being
the ball in the subspace generated by the first $k$
coordinates.

We use $(x',x^{n+1},x'')$ to designate coordinates in
$\bbbr^\nu=\bbbr^n \times\bbbr \times\bbbr^{\nu-n-1}$.
With no loss of generality we may assume that
$$
M\cap\overline{B^\nu}(0,1)=
\overline{B^n}(0,1)\times \{0\}\times\{0\}
\subset
\bbbr^n\times\bbbr\times\bbbr^{\nu-n-1} = \bbbr^\nu.
$$
Let $\gamma$ be as before.
Define $\widetilde{M}$ by replacing
$\overline{B^n}(0,1)\times \{0\}\times\{0\}\subset M$ with the graph
of $\gamma$
$$
\Gamma=
\{(x',x^{n+1},x''):
\ x'\in\overline{B^n}(0,1),\ x^{n+1}=\gamma(x'),\ x''=0\}.
$$
Now we will extend this graph to a graph of a function defined in
$\bbbr^{\nu-1}=\bbbr^n\times\{ 0\}\times\bbbr^{\nu-n-1}$.
Let
$\eta\in C^\infty(\bbbr)$ be such that $0\leq \eta\leq 1$ and
$$
\eta(t) =
\left\{
\begin{array}{ccc}
1   &   \mbox{if $|t|\geq 2$,}\\
0   &   \mbox{if $|t|\leq 1$.}
\end{array}
\right.
$$
Define
$\lambda:\bbbr^n\times\bbbr^{\nu-n-1} \to \bbbr$
by
\begin{equation}
\label{gwiazdka}
\lambda(x',x'') =
\left\{
\begin{array}{ccc}
\gamma(x')\eta\left(\frac{x'}{|x''|}\right)   &   \mbox{if $x''\neq 0$,}\\
\gamma(x')   &   \mbox{if $x''=0$.}
\end{array}
\right.
\end{equation}
It is easily seen that
$\lambda\in C^\infty(\bbbr^n\times\bbbr^{\nu-n-1}\setminus
\{(0,0)\})$,
$\lambda$ is continuous on $\bbbr^{\nu-1}$
and equals zero
outside the ball $B^{\nu-1}(0,2^{-4})$.
Now define the
graph of $\lambda$ over $\overline{B^{\nu-1}}(0,1)$ by
$$
\widetilde{\Gamma} =\{(x',x^{n+1},x'')\in\bbbr^\nu:\
x^{n+1}=\lambda(x',x''),\ (x',x'')\in \overline{B^{\nu-1}}(0,1)\}.
$$
Since $|\lambda|\leq 2^{-6}$ it follows that
$\widetilde{\Gamma}\subset B^\nu(0,2^{-3})$.
Clearly $\widetilde{\Gamma}$ is an extension of $\Gamma$
in the sense that the intersection of
$\widetilde{\Gamma}$ with the linear subspace $x''=0$ coincides with
$\Gamma$.

Observe that
\[
\widetilde{\Phi}(x',x^{n+1},x'') = (x', x^{n+1}-\lambda(x',x''),x'')
\]
is a homeomorphism
$\widetilde{\Phi}:\bbbr^{\nu}\twoheadrightarrow\bbbr^{\nu}$,
$\widetilde{\Phi}(\widetilde{\Gamma})= \overline{B^{\nu-1}}(0,1)$,
$\widetilde{\Phi}(\Gamma)=\overline{B^n}(0,1)$
and
$\widetilde{\Phi}$ is a smooth diffeomorphism off the line
$\{ x'=0,x''=0\}$. However,
$\widetilde{\Phi}$ is {\em not} differentiable in points of the line
$\{ x'=0, x''=0\}$. Moreover $\widetilde{\Phi}$ is not identity
outside any finite ball.
Hence we have to modify $\widetilde{\Phi}$ in order to improve its properties.
Actually we will prove.
\begin{lemma}
\label{king}
There is a $C^\infty$ homeomorphism
$\Phi:\bbbr^{\nu}\twoheadrightarrow\bbbr^{\nu}$ which is diffeomorphism in
$\bbbr^{\nu}\setminus\{ 0\}$,
$\Phi(\widetilde{\Gamma})=\overline{B^{\nu-1}}(0,1)$,
$\Phi(\Gamma)=\overline{B^n}(0,1)$
and $\Phi$ is
identity in $\bbbr^{\nu}\setminus B^{\nu}(0,1)$.
\end{lemma}

It is obvious now that  Lemma~\ref{king} implies
Theorem~\ref{smooth}, so we are left with the proof of the lemma.
Since a version of Lemma~\ref{king} holds for a much larger class
of graphs, and not just for the graph of the function
$\lambda$ constructed above, it makes sense to formulate the lemma
in a more general form and deduce Lemma~\ref{king} as a special
case.
Lemma~\ref{king} is a direct consequence of Lemma~\ref{michalek}
and Lemma~\ref{natusia} below.

Notation $(x',x^{n+1},x'')$ for the coordinates in $\bbbr^\nu$
will no longer be convenient for us, so in the remaining part of
Section~\ref{four} we will denote the coordinates in $\bbbr^\nu$ by
$(x',x^\nu)$, where
$x'=(x^1,\ldots,x^{\nu-1})$.

\begin{lemma}
\label{michalek}
Assume that  a nonnegative function $\lambda$ defined on
$\bbbr^{\nu-1}$ is continuous,
$\supp\lambda\subset B^{\nu-1}(0,2^{-3})$, $\lambda(0)=0$,
and
$\lambda\in C^\infty(\bbbr^{\nu-1}\setminus\{ 0\})$.
Denote the graph of $\lambda$ over $\overline{B^{\nu-1}}(0,1)$ by
$\widetilde{\Gamma}$ i.e.,
$$
\widetilde{\Gamma} =
\{(x',x^\nu)\in\bbbr^\nu:\
x'\in\overline{B^{\nu-1}}(0,1),
x^\nu=\lambda(x')\}.
$$
Assume that the graph of $\lambda$
over $B^{\nu-1}(0,2^{-3})$ is contained in the
ball $B^\nu(0,2^{-3})$. Then there is a homeomorphism
$\Phi^*:\bbbr^\nu\twoheadrightarrow\bbbr^\nu$ of the class
$C^\infty(\bbbr^\nu\setminus\{ 0\})$ with the following properties
\begin{enumerate}
\item $\Phi^*$ is of the form $\Phi^*(x',x^{\nu})=(x',\xi(x',x^{\nu}))$,
\item $\Phi^*(0) = 0$,
\item $\Phi^*$ is a diffeomorphism in $\bbbr^{\nu}\setminus\{ 0\}$,
\item $\Phi^*(\widetilde{\Gamma})=\overline{B^{\nu-1}}(0,1)$,
\item $\Phi^*$ is identity in $\bbbr^{\nu}\setminus B^{\nu}(0,1)$.
\end{enumerate}
\end{lemma}

The mapping $\Phi^*$ constructed in the above lemma has all the
properties needed in Lemma~\ref{king} except for the smoothness at
the origin. The following lemma shows how to modify the mapping
$\Phi^*$ in a way that it will have all the properties needed in
Lemma~\ref{king}.

\begin{lemma}
\label{natusia}
Let $\Lambda:\bbbr^{\nu}\to\bbbr^{\nu}$
be a continuous mapping
such that
$\Lambda(0)=0$ and $\Lambda\in C^{\infty}(\bbbr^{\nu}\setminus\{ 0\})$.
Then there is a
homeomorphism $\Theta:\bbbr^{\nu}\twoheadrightarrow\bbbr^{\nu}$ such that
\begin{enumerate}
\item[(a)] $\Theta\in C^{\infty}(\bbbr^{\nu})$,
\item[(b)] $\Theta$ is a diffeomorphism in $\bbbr^{\nu}\setminus \{0\}$,
\item[(c)] $\Theta(0)=0$ and $\Theta(x)=x$ for all $|x|\geq 1$,
\item[(d)] $\Theta(x)$ is of the form
$$
\Theta(x) =
\left\{
\begin{array}{ccc}
\frac{x}{|x|}\eta(|x|),  & \mbox{if $x\neq 0$,}\\
0 &   \mbox{if $x=0$,}
\end{array}
\right.
$$
for some homeomorphism $\eta:[0,\infty)\twoheadrightarrow
[0,\infty)$and hence
$\Theta(\overline{B^k}(0,1))=\overline{B^k}(0,1)$ for
$k=1,2,\ldots,\nu$.
\item[(e)] $\Theta\circ\Lambda\in C^{\infty}(\bbbr^{\nu})$.
\end{enumerate}
\end{lemma}

If we set $\Lambda=\Phi^*$,
where $\Phi^*$ is defined for $\lambda$ given by (\ref{gwiazdka}),
then $\Phi=\Theta\circ\Phi^*$ has
the properties required in
Lemma~\ref{king}.
Indeed,
it follows from the property (1) of $\Phi^*$ that
$\Phi^*(\Gamma)=\overline{B^n}(0,1)$ and hence
$(\Theta\circ\Phi^*)(\Gamma)=\Theta(\overline{B^n}(0,1))=
\overline{B^n}(0,1)$ by the property (d) of $\Theta$.

\medskip

{\em Proof of Lemma~\ref{michalek}.}
Let $f(x')=4\lambda(x')+|x'|^{2}$. Then
$f\in C^{0}(\bbbr^{\nu-1})\cap C^{\infty}(\bbbr^{\nu-1}\setminus\{0\})$
is strictly positive in $\bbbr^{\nu-1}\setminus\{0\}$,
$f(0)=0$ and $f\geq 4\lambda$ in $\bbbr^{\nu - 1}$.

Fix a function $\vi\in C_0^{\infty}(\bbbr)$ such that
$0<\vi(\tau)<2$ for $0<\tau<1$, $\supp\vi=[0,1]$ and
$\int\vi(\tau)\, d\tau=1$. For $x'\neq 0$ we define
$$
g(x',x^{\nu}) =
\frac{-\lambda(x')}{f(x')}
\vi
\left(
\frac{x^{\nu}+f(x')}{f(x')}
\right)\, .
$$
Then $\int g(x',x^{\nu})\, dx^{\nu}=-\lambda(x')$, the support of
the function $x^{\nu}\mapsto g(x',x^{\nu})$ is contained in
the interval $[-f(x'),0]$, the
function is nonpositive in the interval, but
$g(x',x^{\nu})\geq -1/2$.
For $x'\neq 0$ we also define
\[
h(x',x^{\nu}) =
\frac{\lambda(x')}{f(x')-\lambda(x')}
\vi
\left(
\frac{x^{\nu}-\lambda(x')}{f(x')-\lambda(x')}
\right)\, .
\]
Then the support of the function $x^{\nu}\mapsto h(x',x^{\nu})$ is
contained in the interval
$[\lambda(x'),f(x')]$, the function is nonnegative in the interval and
$\int h(x',x^{\nu})\, dx^{\nu}=\lambda(x')$.

Now we define
$$
\xi(x',x^{\nu}) =
\left\{
\begin{array}{ccc}
\int_{\lambda(x')}^{x^{\nu}}(1+g(x',\tau)+h(x',\tau))\, d\tau
             & \mbox{if $x'\neq 0$,}\\
x^{\nu} &   \mbox{if $x'=0$.}
\end{array}
\right.
$$
The integrand is greater than or equal to $1/2$ and hence the function $\xi$
is strictly increasing in $x^{\nu}$ with
$\partial\xi/\partial x^{\nu}\geq 1/2$. Moreover
$\xi(x',x^{\nu}) = x^{\nu}$ for $|x^{\nu}|\geq f(x')$ and
$\xi(x',x^{\nu})=x^{\nu}-\lambda(x')$ for
$0\leq x^{\nu}\leq\lambda(x')$.

Clearly the mapping $\Phi^*(x',x^{\nu})=(x',\xi(x',x^{\nu}))$ is a homeomorphism
of $\bbbr^{\nu}$ onto $\bbbr^{\nu}$ satisfying claims (1), (2) and (4).
Property (5) follows from the formula for $\xi$ because
$\xi(x',x^{\nu})=x^{\nu}$
whenever $\lambda(x')=0$, and the fact that support of $\lambda$ is concentrated
near $0$.

We still have to prove
that $\Phi^*$ is smooth in $\bbbr^{\nu}\setminus\{ 0\}$ and that the Jacobian of
$\Phi^*$ is positive in $\bbbr^{\nu}\setminus\{ 0\}$.

Smoothness of $\Phi^*$ at $(x'_0,x^{\nu}_0)$, $x'_0\neq 0$ follows from the
definition of $\xi$ and smoothness of $\Phi^*$ at $(0,x^{\nu}_0)$,
$x^{\nu}_0\neq 0$
follows from the
observation that $\Phi^*$ is the identity in a neighborhood of
$(0,x^{\nu}_0)$. The mapping
$\Phi^*$ is a diffeomorphism in $\bbbr^{\nu}\setminus\{0\}$ because direct
computation shows that
$\det D\Phi^*(x',x^{\nu})=\partial\xi/\partial x^{\nu}>0$.
This completes the proof of all the properties of $\Phi^*$.
\hfill $\Box$

\medskip

{\em Proof of Lemma~\ref{natusia}.}
In the proof we will need the following result.
\begin{lemma}
\label{elem}
Let $\{ N_n\}_{n=1}^{\infty}$ be an arbitrary increasing sequence of positive
numbers. Then there is a homeomorphism
\[
\eta:[0,\infty)\twoheadrightarrow [0,\infty)
\]
such that $\eta\in C^{\infty}((0,\infty))$, $\eta'(t)>0$ for $t>0$, $\eta(t)=t$
for $t\geq 1$ and
\begin{equation}
\label{square}
N_n \sup_{0<t\leq 2^{-n}} \sup_{0\leq \ell\leq n}|\eta^{(\ell)}(t)|\to 0
\end{equation}
as $n\to\infty$.
\end{lemma}
{\em Proof.} There is a nonnegative function $\vi\in C_{0}^{\infty}(\bbbr)$
with $\supp\vi=[1/2,2]$ and such that
\begin{equation}
\label{76}
\sum_{n=-\infty}^{\infty}
\vi(2^n t) = 1
\end{equation}
for all $t>0$. Indeed, take an arbitrary function $\widetilde{\vi}\in
C_{0}^{\infty}(\bbbr)$ such that $\supp\widetilde{\vi}=[1/2,2]$, and
$\widetilde{\vi}>0$ in $(1/2,2)$. Then $\vi(t)=\widetilde{\vi}(t)/\alpha(t)$,
where
\[
\alpha(t) = \sum_{n=-\infty}^{\infty} \widetilde{\vi}(2^n t)
\]
satisfies (\ref{76}). This follows from the observation that $\alpha$ is a
strictly positive $C^\infty$ function for all $t>0$ and
$\alpha(2^n t)=\alpha(t)$ for all integers $n$.

Let $\{ a_n\}_{n=2}^{\infty}$ be a sequence of positive numbers less than $1$
such that $a_n\to 0$ as $n\to\infty$. Let also $a_n=1$ for all integers
$n\leq 1$. The function
\[
\psi(t)=\sum_{n=-\infty}^{\infty} a_n \vi(2^n t)
\]
is smooth and strictly positive. Moreover $\psi(t)\to 0$ as $t\to 0+$. Define
\[
\widetilde{\eta}(t) = \int_{0}^{t} \psi(s) \, ds.
\]
Obviously $\widetilde{\eta}:[0,\infty)\twoheadrightarrow [0,\infty)$ is a
homeomorphism,
$\widetilde{\eta}\in C^{\infty}((0,\infty))$ and $\eta'(t)>0$ for $t>0$.
Moreover $\widetilde{\eta}'(t)=1$ for $t\geq 1$. However,
$\widetilde{\eta}(t)<t$ for $t\geq 1$ and hence
$\widetilde{\eta}(t)=t-a$ for all $t\geq 1$ and some constant $a$. Let
$\delta\in C^{\infty}(\bbbr)$ be an arbitrary nondecreasing
function such that $\delta(t)=0$ for $t\leq 1/2$ and $\delta(t)=1$ for $t\geq
3/4$. Now the function $\eta(t)=\widetilde{\eta}(t)+a\delta(t)$ has all the
properties required in the lemma except perhaps for (\ref{square}). We will
show that an appropriate choice of $\{a_n\}_{n=2}^{\infty}$ will guarantee
(\ref{square}) as well.

Since $\delta(t)=0$ for $t\leq 1/2$ we have that $\widetilde{\eta}(t)=\eta(t)$
for small $t$ and hence it suffices to prove (\ref{square}) for
$\widetilde{\eta}$. Let
\[
M_n=\sup_{\ell\leq n-1}\Vert\vi^{(\ell)}\Vert_{\infty}.
\]
We will show that the function $\widetilde{\eta}$ constructed for the sequence
\[
a_n = 2^{-2n^2} N_n^{-1} M_{n}^{-1},
\qquad
n=2,3,4,\ldots
\]
satisfies (\ref{square}).

Since for every $t>0$ no more than two of the functions $\{\vi(2^n
t)\}_{n=-\infty}^{\infty}$ are different from zero we conclude that for
any
$t\leq 2^{-n}$ there is $i\geq n$ such that for
any $1\leq\ell\leq n$
$$
\eta^{(\ell)}(t) =
a_i 2^{i(\ell-1)} \vi^{(\ell-1)}(2^{i}t) +
a_{i+1} 2^{(i+1)(\ell -1)} \vi^{(\ell -1)}(2^{i+1} t).
$$
Hence
$$
|\eta^{(\ell)}(t)| \leq
M_n \left( 2^{-2i^2} 2^{i(\ell-1)} N_n^{-1} M_n^{-1} +
2^{-2(i+1)^{2}}2^{(i+1)(\ell -1)} N_n^{-1} M_n^{-1} \right) \leq
2\cdot 2^{-n^2} N_n^{-1},
$$
and thus
$$
N_n \sup_{0<t\leq 2^{-n}} \sup_{1\leq \ell\leq n} |\eta^{(\ell)}(t)|
\leq 2\cdot 2^{-n^2}\to 0
$$
as $n\to\infty$. It remains to show that
\[
N_n \sup_{0<t\leq 2^{-n}} |\eta(t)|\to 0
\]
as $n\to\infty$. This is easy and we leave it to the reader. The proof of the
lemma is complete.
\hfill $\Box$

\medskip

Now we can complete the proof of Lemma~\ref{natusia}.

We will construct a homeomorphism  $\Theta$ of the form
$$
\Theta(x) =
\left\{
\begin{array}{ccc}
\frac{x}{|x|}\eta(|x|),  & \mbox{if $x\neq 0$,}\\
0 &   \mbox{if $x=0$,}
\end{array}
\right.
$$
where $\eta $ is as in Lemma~\ref{elem} for an appropriate sequence
$\{N_n\}_{n=1}^{\infty}$. The sequence will be chosen later. Obviously
$\Theta:\bbbr^\nu\twoheadrightarrow\bbbr^\nu$ is a homeomorphism such that
$\Theta(0)=0$, $\Theta(x)=x$ for $|x|\geq 1$ and $\Theta$ is a diffeomorphism of
$\bbbr^\nu\setminus\{ 0\}$ of the class $C^\infty$. Observe
that
$$
|\nabla^{m}\Theta(x)| \leq
C(m)\sum_{i=0}^{m} |x|^{-i}|\eta^{(m-i)}(|x|)|,
$$
where $\nabla^m\Theta$ denotes the vector whose
components are all partial derivatives of order $m$. Hence
if $\{ N_n\}_{n=1}^{\infty}$ is a sequence such that
$N_n\geq 2^{n^2}$, then
$$
|\nabla^{m} \Theta(x)|/|x|\to 0
\qquad
\text{as $x\to 0$}
$$
for each positive integer $m$. This, in turn, implies that $\Theta$
has partial
derivatives at $0$ of an arbitrary order and all the derivatives are equal to
$0$. Hence $\Theta\in C^{\infty}(\bbbr^\nu)$.

Obviously $\Theta\circ\Lambda\in C^{\infty}(\bbbr^{\nu}\setminus\{0\})$.
Now it remains to prove that $\Theta\circ\Lambda\in C^{\infty}(\bbbr^{\nu})$ for
some increasing sequence $\{ N_n\}_{n=1}^{\infty}$, such that
$N_n\geq 2^{n^2}$.
To this end it suffices to show that
for some sequence $\{ N_n\}_{n=1}^{\infty}$ we have
\[
|\nabla^{m}(\Theta\circ\Lambda)(x)|/|x|\to 0,
\qquad
\text{as $x\to 0$}
\]
for each nonnegative integer $m$. A multiple application of the Leibniz rule
gives the following variant of the Fa\'a di
Bruno\footnote{ Fa\'a di Bruno was declared a Saint by Pope John
Paul II on September 25, 1988.}
formula
\begin{eqnarray*}
\lefteqn{
|\nabla^{m}(\Theta\circ\Lambda)(x)| \leq
C\left( \sum_{i=1}^{m} |(\nabla^{i}\Theta)(\Lambda(x))|\right) \times} \\
& \times &
\left( \sum |\nabla\Lambda(x)|^{i_1} |\nabla^{2}\Lambda(x)|^{i_2}
\cdot\ldots\cdot|\nabla^{m}\Lambda(x)|^{i_m} \right),
\end{eqnarray*}
where the last sum is taken over all nonnegative integers $i_1$, $i_2$,\ldots,
$i_m$ such that
$$
i_1 + 2 i_2 + 3 i_3 + \ldots + m i_m = m.
$$
Define
$$
M_n = \sup
\left(
\sum
|\nabla\Lambda(x)|^{i_1} |\nabla^{2}\Lambda(x)|^{i_2}
\cdot\ldots\cdot|\nabla^{m}\Lambda(x)|^{i_m} \right) |x|^{-1},
$$
where the supremum is taken over all $x$ satisfying $|x|\leq 1$ and
$2^{-(n+1)} \leq |\Lambda (x)|\leq 2^{-n}$.
For such $x$ we have%
\begin{eqnarray*}
\lefteqn{
|\nabla^{m}(\Theta\circ\Lambda)(x)|/|x| \leq
C M_n \sum_{i=1}^{m}|(\nabla^{i}\Theta)(\Lambda(x))|} \\
& \leq &
C M_n \sum_{i=1}^{m} \left( \sum_{j=0}^{i} |\Lambda(x)|^{-j}
|\eta^{(i-j)}(|\Lambda(x)|)| \right) \leq
C M_n 2^{(n+1)m} m^2 \sup_{\ell\leq m}
|\eta^{(\ell)}(|\Lambda(x)|)|.
\end{eqnarray*}
Hence with a new constant $C$ depending on $m$ for
$n\geq m$ we have
$$
|\nabla^{m}(\Theta\circ\Lambda)(x)|/|x| \leq
C M_n 2^{nm} \sup_{0<t\leq 2^{-n}} \sup_{\ell\leq n}
|\eta^{(\ell)}(t)|.
$$
Now it suffices to take any increasing sequence
$\{ N_n\}_{n=1}^{\infty}$, such that
$N_n\geq M_n 2^{n^2}$.
This completes the proof of Lemma~\ref{natusia} and hence that for
Lemma~\ref{king}. Thus the proof of Theorem~\ref{smooth} is complete.
\hfill $\Box$

\section{Proof of Theorem~\ref{subaru}.}

Let $M=S^n\subset\bbbr^{n+1}$ and let $N=\widetilde{M}$ be the
singular submanifold of $\bbbr^{n+1}$ constructed in the proof of
Theorem~\ref{smooth}.
Let $\vi:S^n=M\to\widetilde{M}=N$ be the $W^{1,n}$ homeomorphism
constructed in the proof of Theorem~\ref{smooth}.
Let $W$ be a closed $n$-dimensional
manifold and let $B\subset W$ be a subset
diffeomorphic to $B^n(0,1)$. Then there is a smooth mapping
$\psi:W\to S^n$ which is one-to-one on $B$ and which maps all of
$W\setminus B$ into one point. Clearly such a mapping has degree
modulo $2$ equal $1$.
Now the same
topological argument based on Lemma~\ref{niema}, as in the proof
of Theorem~\ref{smooth} shows that $\vi\circ\psi\in W^{1,n}(W,N)$
cannot be approximated by Lipschitz mappings $\lip(W,N)$.
\hfill $\Box$

\section{Proof of Theorem~\ref{main}.}

\noindent
{\bf Construction of $X$.}
Let $N$ be a submanifold of $\bbbr^{n+1}$ diffeomorphic to $S^n$ having
$\overline{B^n}(0,1)$ as a subset.
Let $\eta_{k,i}$ be defined as in the proof of Theorem~\ref{smooth} and define
\[
\gamma_{m}(x) = \sum_{k=6}^{m}\sum_{i=1}^{N_k} \eta_{k,i}(x).
\]
Thus the sequence $\{\gamma_m\}$ is a smooth approximation of $\gamma$ defined
by formula (\ref{myszka}) in the $W^{1,n}$ norm.
Let $\Gamma_m$ be the graph of $\gamma_m$ over the unit ball
$\overline{B^n}(0,1)$.

Denote by $\widetilde{N}^m$ the manifold obtained from $N$ by substituting
$\overline{B^n}(0,1)$ with the graph $\Gamma_m$. Since each of the
functions $\gamma_m$ is smooth, $\widetilde{N}^m$
is a smooth manifold diffeomorphic to $S^n$. The
manifolds $\widetilde{N}^m$ converge in some sense to $N^\infty$, the space
obtained from $N$ by
substituting $\overline{B^n}(0,1)$ with $\Gamma$, the graph of
$\gamma$ defined by formula (\ref{myszka}).

Let $B_k$ denote balls defined as in the
definition of the function $\gamma$.
Denote the coordinates in $\bbbr^{n+1}$ by $(x',x^{n+1})$, $x'\in\bbbr^n$,
$x^{n+1}\in\bbbr$.
There is a natural projection
$\pi_m$ of the graph $\Gamma$ onto $\Gamma_m$, namely
$$
\pi_m((x,\gamma(x)) =
\left\{
\begin{array}{ccc}
(x,0)    & \mbox{if $x\in \bigcup_{k=m+1}^{\infty} B_k$,}\\
(x,\gamma(x)) &   \mbox{otherwise.}
\end{array}
\right.
$$
The projection extends in a natural way
to a homeomorphism $\pi_m: N^\infty\twoheadrightarrow \widetilde{N}^m$.

Let $A_m$ be the graph of $\gamma$ over the set
$\bigcup_{k=m+1}^{\infty} B_k$.
Hence $A_m$ is a decreasing sequence of subsets of $\Gamma$ (and hence subsets
of $N^\infty$) with empty intersection.

Denote the orthogonal projection from $\bbbr^{n+1}$ onto $\bbbr^n$
by $\pi$. Observe that
$\pi_m|_{\overline{A}_m}=\pi|_{\overline{A}_m}$.
Hence $\pi_m$ restricted to $\overline{A}_m$ is smooth with the
constant norm of the derivative that does not depend on $m$. Since
$\pi_m$ is identity outside $A_m$, $\pi_m$ extends to a smooth
mapping of a neighbourhood of $N^\infty$ in $\bbbr^{n+1}$ onto
$\widetilde{N}^m$.

\begin{lemma}
\label{dwanascie}
Let $M$ be an $n$-dimensional closed manifold.
If $f\in W^{1,n}(M,N^\infty)$, then
$\pi_m\circ f\in W^{1,n}(M,\widetilde{N}^m)$ and
$\pi_m\circ f\to f$ in $W^{1,n}(M,\bbbr^{n+1})$ as $m\to\infty$.
\end{lemma}
{\em Proof.}
Since $\pi_m$ extends to a smooth mapping we conclude that
$\pi_m\circ f \in W^{1,n}(M,\widetilde{N}^m)$.
Now $\pi_m\circ f\neq f$ on the set $f^{-1}(A_m)$. This
is a decreasing sequence of subsets of $M$ with empty intersection and hence
$\H^n(f^{-1}(A_m))\to 0$ as $m\to\infty$. This, in turn, easily implies that
$\pi_m\circ f\to f$ in $W^{1,n}$ as $m\to\infty$ because
\begin{eqnarray*}
\Vert f-\pi_m\circ f\Vert_{W^{1,n}(M)}
& = &
\left(\int_{f^{-1}(A_m)}| f-\pi_m\circ f|^n \right)^{1/n} \\
& + &
\left(\int_{f^{-1}(A_m)}|D(f-\pi_m\circ f)|^n\right)^{1/n}  \\
& \leq &\left(\H^n(f^{-1}(A_m))\right)^{1/n}
+
\left(\int_{f^{-1}(A_m)}|D(f-\pi\circ f)|^n\right)^{1/n}\\
& \leq &\left(\H^n(f^{-1}(A_m))\right)^{1/n}
+
C \left(\int_{f^{-1}(A_m)}|Df|^{n}\right)^{1/n} \to 0
\end{eqnarray*}
as $m\to\infty$.
We employed here the fact that $|f(x)-(\pi_m\circ f)(x)|\leq 1$
for all $x\in f^{-1}(A_m)$. The proof of the lemma is complete.
\hfill $\Box$

\medskip

$\widetilde{N}^m$ is a subset of $\bbbr^{n+1}$ and hence it is a subset of
$\bbbr^{n+2}=\bbbr^{n+1}\times\bbbr$.
Let $N^m$ be the translation of $\widetilde{N}^m$
by the vector
$\langle 0,\ldots,0,2^{-m}\rangle$ in $\bbbr^{n+2}$. Now we define
$$
\widetilde{X}
= N^\infty \cup \bigcup_{m=6}^{\infty} N^m.
$$
We equip $\widetilde{X}$ with the metric induced from $\bbbr^{n+2}$.
$\widetilde{X}$ is a compact set, but it is not connected.
We define $X$ by connecting all parts of $\widetilde{X}$ by a
curve whose each part is non-rectifiable.

Fix an $n$-dimensional closed manifold $M$ arbitrarily.
We can assume that $M$ is
connected since different components are treated separately.

It is an easy consequence of the absolute continuity of $f$ on lines
\cite[Theorem~2.1.4]{ziemer} that images of almost all lines under
$f$ are rectifiable. Hence Sobolev mappings from $M$ into a space
that contains no rectifiable curves (except constant ones) must be
constant. Accordingly, if $f\in W^{1,n}(M,X)$ is not constant,
then $f\in W^{1,n}(M,\widetilde{X})$.

In the next two steps
we will prove the density of smooth mappings $C^\infty(M,X)$ in $W^{1,n}(M,X)$
and we will construct a global bi-Lipschitz homeomorphism
$\Phi:\bbbr^{n+2}\twoheadrightarrow\bbbr^{n+2}$ such that Lipschitz mappings
$\lip(M,\Phi(X))$ are not
dense in $W^{1,n}(M,\Phi(X))$.

\medskip

\noindent
{\bf Density.}
Constant mappings are smoothm so it suffices to prove the density
of $C^{\infty}(M,\widetilde{X})$ in $W^{1,n}(M,\widetilde{X})$.
We will need the following fact.
\begin{lemma}
For each $f\in W^{1,n}(M,\widetilde{X})$ we can choose a
representative (in the class of functions equal a.e.)
such that $f(M)\subset N^k$ for some $k=6,7,8,\ldots$ or $k=\infty$.
\end{lemma}
{\em Proof.} The lemma is a fairly easy consequence of the absolute continuity
of $f$ on almost all lines.
Indeed, if we would find two subsets $A,B\subset M$
of positive measure such that $f(A)\subset N^{k_1}$ and $f(B)\subset N^{k_2}$,
$k_1\neq k_2$,
then we would find a curve joining $A$ and $B$ such that $f$ is continuous
along that curve. Hence the image of $f$ would contain a curve joining $f(A)$
and $f(B)$ inside $\widetilde{X}$. This is, however, impossible.
\hfill $\Box$

\medskip

Fix $f\in W^{1,n}(M,\widetilde{X})$. We will show that $f$ can be approximated by
mappings from $M$ to $\widetilde{X}$ which are of the class $C^\infty$ as
mappings from $M$ to $\bbbr^{n+2}$.

We know that $f(M)\subset N^k$ for some $k$
finite or $k=\infty$. If $k$ is finite, then there is no problem with the
approximation because $N^k$ is a closed manifold and the smooth mappings
$C^{\infty}(M,N^k)$ are dense in the class $W^{1,n}(M,N^k)$ according
to Schoen and Uhlenbeck's theorem (cf.\ Theorem~\ref{dens}).
Thus we are left with the case in which $f(M)\subset N^\infty$.
Let
\begin{equation}
\label{asterix}
f_m = (\pi_m \circ f) + \langle 0,\ldots,0,2^{-m}\rangle
\end{equation}
i.e. first we compose $f$ with the projection $\pi_m$ from $N^\infty$ onto
$\widetilde{N}^m$ and then we translate the image by the vector
$\langle 0,\ldots,0,2^{-m}\rangle$,
so that the resulting mapping $f_m$ maps $M$ to $N^m$.
It follows from
Lemma~\ref{dwanascie} that $f_m\to f$ in $W^{1,n}$ as $m\to\infty$.
Now it suffices to observe that $f_m\in W^{1,n}(M,N^m)$ can be approximated by
smooth mappings $C^{\infty}(M,N^m)$.
This completes the proof of density of $C^{\infty}(M,\widetilde{X})$ in
$W^{1,n}(M,\widetilde{X})$.

\vspace{3mm}

\noindent
{\bf Construction of $\Phi$.}
We can assume that $[1,2]^n$ is a subset of $N$. Since
$\widetilde{N}^m$ was obtained from $N$ by modifying it on
$B^n(0,1)$ only, $[1,2]^n$ is a subset of $\widetilde{N}^m$ as
well.  Thus
$$
I_m=[1,2]^n\times \{ 0\} \times \{ 2^{-m}\}
$$
is a subset of $N^m$. Moreover the intersection of some
$(n+2)$-dimensional box containing $I_m$ with
$\widetilde{X}$ equals $I_m$.
More precisely we can assume that
$$
\widetilde{X}\cap K_m=I_m,
$$
where
$$
K_m=
[1,2]^n\times
[-2^{-(m+2)},2^{-(m+2)}]\times
[2^{-m}-2^{-(m+2)},2^{-m}+2^{-(m+2)}].
$$
Now for each $m$ we will construct an $L$-bi-Lipschitz homeomorphism
$\Phi_m:K_m\twoheadrightarrow K_m$ which is identity on the boundary
$\partial K_m$, and with the constant $L$ that does not depend on $m$.
It follows that the mapping
$$
\Phi(x)=
\left\{
\begin{array}{ccc}
\Phi_m(x)    & \mbox{if $x\in K_m$, $m=6,7,8,\ldots$}\\
x &   \mbox{otherwise.}
\end{array}
\right.
$$
is bi-Lipschitz. Indeed, it is enough to observe that the norm of
the derivative is bounded from below by a positive constant and
bounded from above.

The set $\Phi_m(I_m)$ will be a subset of $K_m$
with ``wrinkles''. It will be a graph of a Lipschitz function
$\vi_m$ constructed later.
Thus the mapping $\Phi$ transforms each $N^m$ into a manifold with
``wrinkles'' added to its surface.
Note that $\Phi(N^\infty)=N^\infty$.

Let $M$ be a closed $n$-dimensional manifold.
According to Theorem~\ref{subaru}
there is a mapping $f\in W^{1,n}(M,N^\infty)$ which cannot be
approximated by Lipschtz mappings $\lip(M,N^\infty)$.
Since $\Phi(N^\infty)=N^\infty$, we have
$f\in W^{1,n}(M,\Phi(X))$.
Our aim is to construct a
bi-Lipschitz mapping $\Phi:\bbbr^{n+2}\twoheadrightarrow\bbbr^{n+2}$
of the form described above which has the property that
there is $\eps>0$ such that every Lipschitz mapping
$f_m:M\to \Phi(N^m)$ satisfies
$\Vert f-f_m\Vert_{1,n}>\eps$.
This will readily imply that the mapping
$f\in W^{1,n}(M,\Phi(X))$ cannot be approximated by
Lipschitz mappings $\lip(M,\Phi(X))$.
Note that $f$ can be approximated by Lipschitz mappings into
$N^m$ (see (\ref{asterix})), however, adding ``wrinkles'' to $N^m$ by applying $\Phi$
will have the effect that mappings into $\Phi(N^m)$ stay away from $f$.

Now we will construct Lipschitz functions $\vi_m$ whose graphs will be the sets
$\Phi(I_m)$.

Let
$E_m=\{x\in K_m:\, x_{n+1}=0\}$
be an $(n+1)$-dimensional box contained in $K_m$
and let
$E_m^i=\{ x\in K_m:\, \mbox{$x_{n+1}=0$ and $x_n=1+i2^{-(m+10)}$}\}$
for $i=0,1,2,\ldots,2^{m+10}$.
The set
$$
F_m=\partial E_m\cup\bigcup_{i=0}^{2^{m+10}} E_m^i
$$
is an $n$-dimensional subset of $E_m$.
Let $\vi_m:E_m\to\bbbr$ be the function defined by the formula
$$
\vi_m(x)=\dist (x,F_m)
\qquad
\mbox{for $x\in E_m$}.
$$
Using notation $(x',x^{n+1},x^{n+2})$, $x'\in\bbbr^n$, for the coordinates
in $\bbbr^{n+2}$, we define the graph of $\vi_m$ by
$$
G_m=\{(x',\vi_m(x',0,x^{n+2}),x^{n+2}):\,
(x',0,x^{n+2})\in E_m\}.
$$
$G_m$ is a compact subset of $K_m$.
The Lipschitz constant of $\vi_m$ is $1$ as it is the distance function.

\begin{lemma}
\label{joasia}
There is a constant $L>1$ depending on $n$ only such that for
every integer $m=6,7,8,\ldots$ there is an $L$-bi-Lipschitz homeomorphism
$\Phi_m:K_m\twoheadrightarrow K_m$
which is identity on the boundary and which satisfies
$\Phi_m(E_m)=G_m$.
\end{lemma}
{\em Proof.}
Let
$E_m\times\bbbr = \{(x',x^{n+1},x^{n+2}):\,
(x',0,x^{n+2})\in E_m\}$.
With this notation
$K_m=E_m\times[-2^{-(m+2)},2^{-(m+2)}]\subset E_m\times\bbbr$.
The mapping
$\widetilde{\Phi}_m:E_m\times\bbbr\twoheadrightarrow E_m\times\bbbr$
given by
$\widetilde{\Phi}_m(x',x^{n+1},x^{n+2})=
(x',x^{n+1}+\vi_m(x',0,x^{n+2}),x^{n+2})$
is $L$-bi-Lipschitz for some constant $L$ that does not depend on $m$.
Indeed, the inverse mapping is
$\widetilde{\Phi}_m^{-1}(x',x^{n+1},x^{n+2})=
(x',x^{n+1}-\vi_m(x',0,x^{n+2}),x^{n+2})$
and both $\widetilde{\Phi}_m$ and $\widetilde{\Phi}_m^{-1}$ have partial derivatives
bounded by
$1$, because the Lipschitz constant of $\vi_m$ is $1$.
The mapping $\widetilde{\Phi}_m$ is identity on the part of the boundary of $K_m$
(because $\vi_m$ is zero on $\partial E_m$), however, it is not
identity on the entire boundary of $K_m$.
Fortunately, one can easily use an argument
similar to that employed in the proof of
Lemma~\ref{michalek} and construct a mapping $\Phi_m$ of the form
$\Phi_m(x',x^{n+1},x^{n+2})=(x',\xi_m(x',x^{n+1},x^{n+2}),x^{n+2})$
being a modification of $\widetilde{\Phi}_m$ and satisfying the
claim of Lemma~\ref{joasia}.
Indeed, one can write the fomula for
$\xi_m$ explicitly making it linear as a function of $x^{n+1}$
on each of the segments
$[-2^{-(m+2)},0]$ and  $[0,2^{-(m+2)}]$.
We leave easy details to the
reader.
\hfill $\Box$

\vspace{3mm}

\noindent
{\bf Lack of density.}
Now we can complete the proof of Theorem~\ref{main} by showing that
$f\in W^{1,n}(M,\Phi(X))$ cannot be approximated by mappings $\lip(M,\Phi(X))$.
As we have already observed it suffices to prove that there is no sequence
$f_{m_k}\in\lip(M,\Phi(N^{m_k}))$, such that
$\Vert f-f_{m_k}\Vert_{1,n}\to 0$ as $k\to\infty$.
Let us briefly recall the construction of the mapping
$f$ (cf.\ the proof of Theorem~\ref{subaru};
in that proof $W$ plays the role of $M$).
$N^\infty$ is homeomorphic to $S^n$ and $N^\infty\setminus \{ 0\}$ is
diffeomorphic to $S^n$ with one point removed.
We choose a subset $\overline{B}$ of $M$ diffeomorphic to a closed Euclidean ball.
The mapping $f$ maps $M\setminus B$ into one point
$0\neq x^*\in N^\infty$ and
$f|_B:B\to N^\infty\setminus\{ x^*\}$ is a homeomorphism. Moreover
$f|_{B\setminus f^{-1}(0)}:B\setminus f^{-1}(0)\to N^\infty\setminus \{ x^*,0\}$
is a diffeomorphism.
We can assume that $[1,2]^n\subset N^\infty\setminus \{ x^*,0\}$
as the point
$0\neq x^*\in N^\infty$ can be chosen
arbitrarily.
That means, applying a suitable change of variables in $M$,
we can assume that $[1,2]^n\subset M$ and that $f:M\to N^\infty$ satisfies
$f(x)=x$ for $x\in [1,2]^n$.
Suppose that
$f_{m_k}\in\lip(M,\Phi(N^{m_k}))$,
$\Vert f-f_{m_k}\Vert_{1,n}\to 0$ as $k\to\infty$.
Let us denote the coordinates in $[1,2]^n$ by
$$
x=(\hat{x},x^n)\in [1,2]^{n-1}\times [1,2]=[1,2]^n.
$$
Fubini's theorem implies that for a suitable subsequence of $f_{m_k}$
(also denoted by $f_{m_k}$) and for almost all $\hat{x}\in [1,2]^{n-1}$, the sequence
of functions of one variable
\begin{equation}
\label{B421}
[1,2]\ni x^n\mapsto f_{m_k}(\hat{x},x^n)
\end{equation}
converges to the function on one variable
\begin{equation}
\label{B422}
[1,2]\ni x^n\mapsto f(\hat{x},x^n)
\end{equation}
in $W^{1,n}([1,2])$. In particular it implies that the functions converge uniformly
and that the lengths of the curves (\ref{B421}) converge to the length of the curve
(\ref{B422}) as $k\to\infty$.
The length of the curve (\ref{B422}) equals
$1$, because $f$ is identity on $[1,2]^n$.
However, for $k$ sufficiently large, the length of the curve (\ref{B421}) is
at least $\sqrt{2}$, because it goes across the ``wrinkles'' $\Phi(I_{m_k})$.
Indeed, the projection of the curve (\ref{B421}) onto the coordinate plane
generated by coordinates $(x^n,x^{n+1})$ contains the projection of the set
$\Phi(I_{m_k})$ on that plane. This projection is, however, isometric to the curve
$[1,2]\ni t\mapsto (t,\dist (t,S))$,
where $S=\{1+i2^{-(m+10)}:\, i=0,1,2,\ldots,2^{m+10}\}$
and the length of that curve equals $\sqrt{2}$.
Since the orthogonal projection does not increase the length, we conclude that the
length of the curve (\ref{B421}) is at least $\sqrt{2}$.
This is a contradiction with the convergence of the length.
The proof of Theorem~\ref{main} is complete.
\hfill $\Box$

\medskip

{\em Remark.}
The homeomorphisms $\Phi_m:K_m\twoheadrightarrow K_m$
are not smooth. However, one can easily modify the construction in a way that
the they are diffeomorphism. Then the mapping
$\Phi:\bbbr^{n+2}\twoheadrightarrow\bbbr^{n+2}$ is a diffeomorphism
in $\bbbr^{n+2}\setminus [1,2]^{n}$ i.e. it is a diffeomorphism outside a
square of dimension $n$.

\end{document}